\documentclass{article}
\usepackage{epsfig}
\usepackage{latexsym}
\usepackage{amsmath}
\usepackage{amsthm}
\usepackage{amssymb}
\usepackage{graphicx}
\usepackage{thesis}
\usepackage{eufrak}
\begin{document}
\title{A TQFT of Intersection Numbers on Moduli Spaces of Admissible Covers}
\author{Renzo Cavalieri}

\maketitle
\begin{abstract}
We construct a two-level weighted TQFT whose structure coefficents are equivariant intersection numbers on moduli spaces of admissible covers. Such a structure is parallel (and strictly related) to the local Gromov-Witten theory of curves in \cite{bp:tlgwtoc}. We compute explicitly the theory using techniques of localization on moduli spaces of admissible covers of a parametrized $\proj$. The Frobenius Algebras we obtain are one parameter deformations of the class algebra of the symmetric group $S_d$. In certain special cases we are able to produce explicit closed formulas for such deformations in terms of the representation theory of $S_d$. 
\end{abstract}
\section*{Introduction}
This paper studies a large class of (equivariant) intersection 
numbers on moduli spaces of admissible covers. For a smooth algebraic curve $X$, ramified covers of a given degree  of $X$ by smooth curves of a given genus are parametrized by moduli spaces called Hurwitz schemes. A  smooth compactification of a Hurwitz scheme can be obtained by allowing suitable degenerations, called admissible covers.

Moduli spaces of admissible covers were introduced originally by Harris
and Mumford in
\cite{hm:kd}.  Intersection theory on these spaces was for a long time
extremely hard and mysterious, mostly because they are in general not
normal, even if the normalization is always smooth.  Only recently in
\cite{acv:ac}, Abramovich, Corti and Vistoli exhibit this normalization
as the stack of balanced stable maps of degree 0 from twisted curves to
the classifying stack $ \mathcal{B}S_d $.  This way they attain both the
smoothness of the stack and a nice moduli-theoretic interpretation of
it.  We abuse notation and refer to the Abramovich-Corti-Vistoli (ACV)
spaces as admissible covers.

At about  the same time, Ionel developed a parallel theory in the
symplectic category (\cite{i:trrih}) and used push-pull techniques on admissible covers to produce new relations in the tautological ring of $M_{g,n}$ (see also \cite{i:rittr}). 
   
In \cite{gv:rvl}, admissible cover loci within the boundary of moduli spaces of stable maps play a key role in establishing the result that the degree $3g-3$ part of the tautological ring of $\overline{M}_g$ has dimension $1$, providing further evidence for a conjecture by Faber, stating that $R(\overline{M}_g)$ is a Gorenstein algebra with socle in degree $3g-3$.

Recently, in \cite{bgp:crc}, Bryan, Graber and Pandharipande show that the orbifold Gromov-Witten potential of a  Gorenstein orbifold can be computed in terms of intersection theory on moduli spaces of admissible covers. With a subtle use of WDVV techniques, they are able to explicitly compute  the Gromov-Witten potential for the orbifold $[\mathbb{C}^2/\mathbb{Z}_3]$. Such a computation provides evidence for the crepant resolution conjecture (\cite{bg:prep}).

We give a few basic definitions and a working description of moduli spaces of admissible covers in section \ref{admcov}.

\vspace{0.5cm}

For all choices of:
\begin{itemize}
        \item an $r$-pointed curve $(X,p_1,\ldots,p_n)$;
        \item a rank two vector bundle $N=L_1\oplus L_2$ on $X$, endowed with a natural 
$\mathbb{C}^\ast \times \mathbb{C}^\ast$ action (page \pageref{aci});
        \item  a vector of partitions $\underline{\eta}=(\eta_1,\ldots,\eta_n)$,
\end{itemize}
we describe the invariants
$$
A_d^h(N):=\int_{ \overline{Adm}_{h\stackrel{d}{\rightarrow}X,
(\eta_1 p_1,\ldots,\eta_rp_r)}}
 \hspace{-1.5cm} e^{eq}(-R^\bullet\pi_\ast f^\ast(L_1\oplus L_2)).
$$

Motivation for studying these invariants is twofold: they are natural and interesting intersection numbers on their own, that give rise to a beautiful structure. In the context of Gromov-Witten theory, invariants of this form are known as ``local" invariants: roughly speaking, they represent the contribution to the Gromov-Witten invariants of a threefold given by rigid curves.

\vspace{0.5cm}
\textbf{Theorem \ref{tqftteor}}  (page \pageref{tqftteor}) \textbf{:} \textit{The invariants  $A_d^h(N)$ can be organized to be the structure coefficients of a
\textbf{$2-$level, semisimple, weighted TQFT $\mathcal{U}$}.}
\vspace{0.5cm}

Section \ref{tqft} is dedicated to presenting these structures to the unfamiliar reader, while in section
\ref{u} the specific TQFT $\mathcal{U}$ is constructed.

The generators for the TQFT are explicitly computed in section \ref{compute}. The techniques involved are basic dimension counts, reduction to classical intersection theory on moduli spaces of curves, and Atiyah-Bott localization on moduli spaces of admissible covers of a parametrized $\proj$.

An interesting feature of this theory is that the degree $0$ part is constructed from Hurwitz numbers.
The embedded (see page \pageref{emb}) Frobenius Algebras induced on the Hilbert space by $\mathcal{U}$
are one parameter deformations of the class algebra of the symmetric group, whose
TQFT-theoretic description in terms of Hurwitz numbers was studied in the '90s in \cite{dw:tgtagc} and \cite{fq:cst}.
An explicit description of such deformations is in general quite complicated. By specializing  
 to the anti-diagonal action
of $\mathbb{C}^\ast$ inside $\mathbb{C}^\ast\times\mathbb{C}^\ast$, it is possible
to diagonalize the theory: closed formulas for our invariants and for the deformation are   
described in terms of the representation theory of the symmetric group $S_d$ (\textbf{Theorem~\ref{antid}}, page \pageref{antid}). This is carried out in section \ref{special}.

This work is closely connected to and follows recent work of Jim Bryan and Rahul 
Pandharipande (\cite{bp:tqft},\cite{bp:tlgwtoc}), describing the local Gromov-Witten theory of curves.

There, analogous intersection numbers on moduli spaces of (relative) stable maps are organized
in a TQFT, that we  denote $\mathcal{BP}$.
 \textbf{Theorem} \ref{equals} (page \pageref{equals}) shows that the two theories coincide 
in level $(0,0)$. In all other levels, $\mathcal{U}$
is a normalization of $\mathcal{BP}$ via appropriate powers of a universal 
generating function factor, which should be understood as the contribution of maps containing contracting
components to the Gromov-Witten invariants.

This fact, the most technical result in this paper, is established
by computing the genus $0$, one-pointed invariants via localization, together with the use
of some beautiful Hodge integral computations by Faber-Pandharipande (\cite{fp:hiagwt}) and 
Ekedahl-Lando-Shapiro-Vainshtein (\cite{elsv:hnaiomsoc} and \cite{vg:hnavl}). The explicit result is:

\vspace{0.5cm}
\textbf{Theorem \ref{CYcap}} (page \pageref{CYcap}) \textbf{:}\textit{the coefficients for the one-pointed invariants of  $\mathcal{U}$  in level $(0,-1)$ are given by the following generating functions:
$$
A_d(0|0,-1)_{\eta}=(-1)^{d-\ell(\eta)}\frac{\left(2\sin\left(\frac{u}{2}\right)\right)^d}{(s_1)^{\ell(\eta)} \mathfrak{z}{(\eta)}\prod 2\sin\left(\frac{\eta_iu}{2}\right)}.
$$
}
\vspace{0.5cm}

A direct check in the one-pointed case, together with the semisimplicity of both theories, yield the following corollary:

\vspace{0.5cm}
\textbf{Corollary 1:} \textit {the coefficients of the theories $\mathcal{U}$ and $\mathcal{BP}$  are related by:}
$$ A_d(g\mid k_1 , k_2)_{\underline{\eta}}=(d!)^{k_1+k_2}s_1^{dk_2}s_2^{dk_1}\mathcal{BP}_d(g,\mid k_1, k_2)_{\underline{\eta}}\mathcal{BP}_d(0\mid0,-1)_{(1,\ldots,1)}^{k_1+k_2}$$
 
This close proximity to Gromov-Witten theory reinforces our interest in 
moduli spaces of admissible covers, as it anticipates the possibility of a fertile exchange of information between the two contexts. In particular, embedded in the theory $\mathcal{U}^\circ$ (the circle superscript indicates we are restricting our attention to connected covers) we rediscover the classical result:

\begin{description}
	\item[\textit{Aspinwall-Morrison formula:}]
	$$
	\int_{[\overline{M}_{0,0}(\proj,d)]} \hspace{-1cm}R^1\pi_\ast f^\ast(\mathcal{O}(-1)\oplus\mathcal{O}(-1))=
	\left(\frac{A_d^{\circ,0}(0 \mid -1,-1)}{u^{2d-2}}\right)_{\mid u=0} = \frac{1}{d^3}
	$$ 
\end{description}

The technique of Atyah-Bott localization suits very well the spaces of admissible 
covers of a parametrized $\proj$; the fact that these spaces are smooth
(as DM stacks) requires no need for a  virtual fundamental class in order to do
intersection theory on them. The modularity of the boundary-fixed loci
naturally produces topological recursions that live completely within the
realm of admissible covers.

\subsection*{Acknowledgements}
This paper owes a lot to Jim Bryan and Rahul Pandharipande, whose work
provided amazing motivation and a useful roadmap. In particular I thank
Jim Bryan for suggesting to me to look at moduli spaces of admissible covers,
and for graciously answering all my pestering emails.
I am grateful to my advisor, Aaron Bertram,
 for his constant support, motivation, and expert guidance.
 I also thank  Alastair Craw, Bill Fulton, Y.P. Lee and Ravi Vakil for extremely useful
conversations and feedback.

\section{Admissible Covers}\label{admcov}
Moduli spaces of admissible covers are a ``natural'' compactification of
the Hurwitz schemes, parametrizing ramified covers of smooth Riemann
Surfaces.  The fundamental idea is that, in order to understand limit
covers, we allow the base curve to degenerate together with the cover.
Branch points are not allowed to ``come together''; as two or more
branch points tend to collide, a new component of the base curve sprouts
from the point of collision, and the points transfer onto it.
Similarly, upstairs the cover splits into a nodal cover.

More formally:  let $ (X,p _1,\cdots,p_r) $ be an $ r $-pointed
nodal curve of genus $ g $.
\begin{defi}
    An \textbf{admissible cover} $ \pi:E\longrightarrow X $ of degree d
    is a finite morphism satisfying the following:
    \begin{enumerate}
        \item
            E is a nodal curve.
        \item
            Every node of E maps to a node of X.
        \item
            The restriction of $ \pi:E\longrightarrow X $ to $ X\setminus
            (p_1,\cdots,p_r) $ is \'{e}tale of constant degree d.
        \item
            Nodes can be smoothed. This means: given an admissible cover
$\pi:E\rightarrow X$, and a node of $E$, we can find a family of admissible 
covers $\pi':E'\rightarrow X'$ such that:
\begin{itemize}
\item $\pi:E\rightarrow X$ is the central fiber of the family;
\item locally in analytic coordinates, $X', E'$ and $ \pi'
            $ are described as follows, for     some positive integer $n$ not larger than $d$:
            $$
                \begin{array}{cl}
                    E:  & e_1e_2=a, \\
                    X:  & x_1x_2=a^n ,\\
                    \pi:  & x_1= e_1^n,\ x_2= e_2^n.
                \end{array}
            $$
\end{itemize}

    \end{enumerate}
\end{defi}

We recall here the notation we use in this paper, and refer the reader to \cite{r:thesis}
for a more extensive discussion.

Let $(X, p_1,\ldots,p_r)$ be as before, and
$\eta_1,\ldots,\eta_r$ be partitions\footnote{
See the appendix for partition notation.}
 of the fixed integer $d$. 
We denote by
$$
\overline{Adm}_{h\stackrel{d}{\rightarrow}X,(\mu_1p_1,\cdots,\mu_rp_r)}
$$
the stack of \textbf{possibly disconnected}, degree $d$ admissible covers of the curve $X$ 
by curves of genus $h$, such that:
\begin{itemize}
\item the ramification profile over the base point $p_i$ is described by the partition
$\eta_i$;
\item all other ramification is required to be simple and is not marked.
\end{itemize}
The following variations are also used:
\begin{description}
\item[connected admissible covers:] we add the superscript ``$\circ$" to restrict our attention
to admissible covers by connected curves;
\item[admissible covers of a genus g curve:] we denote by
$$
\overline{Adm}_{h\stackrel{d}{\rightarrow}g}
$$
the stack of admissible covers of a curve of genus $g$. This means that also
the base curve is allowed to vary in families.
\item[admissible covers of a parametrized $\proj$:] when we intend to fix
a parametrization on the base $\proj$, we write
$$
\overline{Adm}_{h\stackrel{d}{\rightarrow}\proj}.
$$
\end{description}
Moduli spaces of admissible covers admit forgetful maps to (quotients of)
the configuration spaces of points on the base curve (resp. $\overline{M}_{g,n}$ for admissible covers of a genus $g$ curve)- recording the information 
of the branch points that are free to move. 
Tautological  $\psi$ classes on admissible covers
are defined by pulling back the $\psi$ classes downstairs  via these maps.

There is also a natural map from a moduli space of admissible covers of genus $h$ to the corresponding moduli space of curve $\overline{M}_h$, obtained by forgetting the cover map and only remembering the source curve.
Tautological $\lambda$ classes on admissible covers are defined by pulling back $\lambda$ classes (the Chern classes of the Hodge bundle on the moduli space of stable curves) via these maps.

\subsection{Admissible Covers of a Nodal Curve}\label{nodal}

Admissible covers of a nodal curve can be  described combinatorially in terms of 
admissible covers of the irreducible components of the curve. This is extremely 
useful because it opens the way to the use of degeneration techniques and induction.
Crucial to this work are the following identities (\cite{l:adffgwi}), that take place in the Chow ring
with rational coefficients.
\begin{description}
        \item[Reducible nodal curve:] let
$$
    X= X_1\bigcup_{x_1=x_2} X_2
$$
be a nodal curve of genus $ g $, obtained by attaching at a point two
irreducible curves of genus $ g_1 $ and $ g_2 $. Then:

\begin{eqnarray}
        \displaystyle{[\overline{Adm}_{h\stackrel{d}{\rightarrow}X}]=
        \sum_{\eta,h_1,h_2}\mathfrak{z}(\eta)[\overline{Adm}_{h_1\stackrel{d}{\rightarrow}X_1
        ,( {\eta})}] \times [ \overline{Adm}_{h_2\stackrel{d}{\rightarrow}X_2,
        ({\eta})}],} 
\label{rnodal}
 \end{eqnarray}
where:
\begin{itemize}
        \item $\eta = ((\eta^1)^{m_1},\ldots,(\eta^k)^{m_k})$ runs over all partitions
 of $d$;
        \item we define the combinatorial factor: 
\begin{eqnarray}
    \mathfrak{z}(\eta):=\prod m_i!(\eta^i)^{m_i}
    \label{zeta}.
\end{eqnarray}
In particular, $\mathfrak{z}(\eta)$ is the order of the centralizer in $S_d$ of any group element in the conjugacy class of $\eta$;
        \item $h_1+h_2+\ell(\eta)-1=h.$
\end{itemize}

\textbf{Note:} if we are dealing with  admissible cover spaces with also a
prescribed vector of ramification conditions $ \underline{\mu} $, analogous
formulas hold;
the $ \mu_i $'s  need to be distributed on the two twigs $
X_1 $ and $ X_2 $ in all possible ways.

        \item[Irreducible nodal curve:] let
$$
    X= X'/\{x_1=x_2\}
$$
be a nodal curve of genus $ g $, obtained by gluing two distinct points
of an irreducible curve $X'$ of genus $ g-1 $. As an element in the Chow ring with rational coefficients, we can then
express:
\begin{eqnarray}
        \displaystyle{[\overline{Adm}_{h\stackrel{d}{\rightarrow}X}]=
        \sum_{\eta}\mathfrak{z}(\eta)[\overline{Adm}_{h'\stackrel{d}{\rightarrow}X'
        ,(\eta, {\eta})}]},
\label{irrnodal}
\end{eqnarray}
where the sum is over all partitions $\eta$ of $d$, and $h'$ is determined by the relation:
$$
h'+\ell(\eta) = h.
$$

\end{description}

\section{Topological Quantum Field Theories}\label{tqft}

An excellent and elementary reference for two-dimensional topological quantum field theories in mathematics is \cite{j:tqft}. 

\begin{defi}
A \textit{(1+1)
dimensional topological quantum field theory} is a functor of tensor
categories:
$$
    \mathcal{T}:\textbf{2Cob}\longrightarrow \textbf{ Free Rmod}.
$$
\end{defi}
On the right hand side, we have the  category of 
free modules over a commutative ring $R$.
 Let us  describe the
category $ \mathbf{2Cob} $:
\begin{description}
    \item[objects:]
        objects are one-dimensional oriented closed manifolds, i.e.,
        finite disjoint unions of oriented circles.
    \item[morphisms:]
        morphisms are (equivalence classes of) oriented cobordisms
        between two objects.  We can think of them as oriented
        topological surfaces with oriented boundary components.
    \item[composition:]
        we compose two morphisms by simply concatenating them;
        equivalently, we glue negatively oriented boundary components of
        one surface to positively oriented boundary components of the
        other.
    \item[tensor structure:]
        the tensor operation is disjoint union.
\end{description}

The free module $ H:=\mathcal{T}
(S^1) $  is called the \textbf{Hilbert space}
of the TQFT\@.

All topological surfaces can be decomposed into discs, annuli, and pairs
of pants.  Therefore, the structure of a TQFT is completely determined
if it is described on these basic building blocks. 

\subsection{Tensor Notation} 
\begin{figure}[htbp]
    \begin{center}
        \includegraphics[width=12cm]{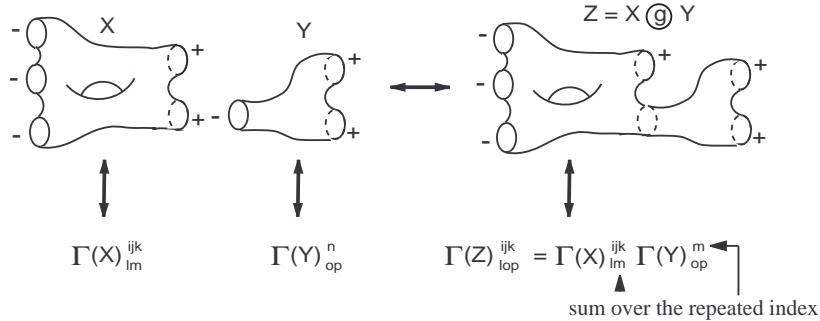}
        \caption{Gluing in tensor notation.}
        \label{glue}
    \end{center}
\end{figure}

It is convenient, for explicit
computations, to familiarize ourselves with tensor notation for TQFT's.
Let us  explicitly choose a basis e$ _1,\dots, $e$ _r $ for the Hilbert
space $H$, and let us denote the dual basis by e$ ^1,\dots, $e$ ^r $. Let $
W_m^n(g) $ be a genus $ g $ cobordism from $ m $ to $ n $ circles.  Then the
map
$$
    \mathcal{T}(W_m^n(g)):H^{\otimes m}\rightarrow H^{\otimes n}
$$
can be thought of as a vector in $ (H^\ast)^{\otimes m} \otimes H^{\otimes
n} $. We  denote by
$$
    \Gamma(W_m^n(g))_{i_1,\dots,i_m}^{j_1,\dots,j_n}
$$
the coefficient of $ \mathcal{T}(W_m^n(g)) $ in the direction of the basis
element e$ ^{i_1}\otimes\dots\otimes $e$ ^{i_m}\otimes $e$ _{j_1}\otimes\dots\otimes
$e$ _{j_n} $.  That is,
$$
    \mathcal{T}(W_m^n(g))=\sum\Gamma(W_m^n(g))_{i_1,\dots,i_m}^{j_1,\dots,j_n}\mbox{e}^
    {i_1}\otimes\dots\otimes\mbox{e}^{i_m}\otimes\mbox{e}_{j_1}\otimes\dots\otimes
    \mbox{ e}_{j_n},
$$
as depicted in Figure~\ref{glue}.

\subsection{Frobenius Algebras} 

A TQFT gives the Hilbert space $ H $ the
structure of a commutative Frobenius algebra.  This means it defines an
associative and commutative multiplication `` $ \cdot\ $'' and an inner
product (also called the metric of the TQFT) ``$ <,> $'' on $ H $ such that
\begin{eqnarray}
    <h_1\cdot h_2,h_3>= <h_1,h_2\cdot h_3>
\end{eqnarray}
holds for all $ h_1 $,$ h_2 $,$ h_3 $ in the Hilbert space $ H $. It is easy
to see how the structure is induced:  multiplication is the map
associated to the $ (-,-,+) $-pair of pants, the inner product is the
scalar map associated to the $ (-,-) $-annulus.  As a consequence, we
see immediately that the cap with positively oriented boundary
corresponds to the unit vector for the multiplication map just defined, whereas the $(-)$-cap corresponds to the counit operator in the Frobenius Algebra.
\subsection{Semisimple TQFT's}\label{semisimple}

\begin{defi}
    A TQFT $ \mathcal{T} $ is semisimple if the Frobenius algebra
    induced on the Hilbert space $ H $ is semisimple.  That is, if there
    is an orthonormal basis $ e_1,\dots $e$ _r $ for $ H $ such that
    $$
        \mbox{e}_i \cdot \mbox{e}_j =\delta_{ij}\mbox{e}_i.
    $$
\end{defi}

An equivalent point of view is to say that $ \mathcal{T} $ is a direct
sum
$$
    \mathcal{T}= \mathcal{T}_1\oplus\dots\oplus\mathcal{T}_r,
$$
where all $ \mathcal{T}_i $'s are TQFT with Hilbert
space equal to the ground ring.

Denote by e$ _1,\dots,  $e$ _r $ a semisimple basis for $ H $.  We can
also think of e$ _i $ being the identity vector for the space $ H_i $.
Let e$ ^1,\dots,  $e$ ^r $ be the dual basis.
Then semisimplicity is equivalent to asking all non-diagonal
coefficients to vanish:
$$
    \Gamma_{i_1,\dots,i_n}^{j_1,\dots,j_m}(W_m^n(g))=0,
$$
unless $ i_1= i_2=\dots=i_n=j_1 = \dots = j_m $.

There are now $ r $ universal constants $ \lambda_1,\dots, \lambda_r $
that govern the structure of the TQFT\@.  They can be defined in many
equivalent ways.  Here are two equivalent descriptions that we will be using later
on:
\begin{enumerate}
    \item
        $ 1/\lambda_i $ is the image of the basis vector e$ _i $ via the
        counit operator.
    \item
        $ \lambda_i $ is the $ i $-th eigenvalue of the genus adding
        operator (this is the linear map associated to the torus with a negative and a positive puncture, represented in Figure \ref{ga}).
\end{enumerate}

Now the following structure theorem holds:
\begin{fact}
    Let $ \mathcal{T} $ be a semisimple TQFT, and all notation  as
    above. Denote by $ W_m^n(g) $ a genus $ g $ surface with $ m $ input
    and $ n $ output holes. Then:
    $$
        \mathcal{T}(W_m^n(g))= \sum_{i=1}^r \lambda_i^{g+n-1}\underbrace{\mbox
        {e}^i\otimes\dots\otimes\mbox{e}^i}_{m\ %
        \ %
        times}\underbrace{\otimes\mbox{e}_i\otimes\dots\otimes\mbox{e}_i}_
        {n \ \ {times}}.
    $$
    In particular:
    \begin{eqnarray}
        \mathcal{T}(W_0^0(g))= \sum_{i=1}^r\lambda_i^{g-1}.
    \end{eqnarray}
\end{fact}

\subsection{The TQFT of Hurwitz Numbers}\label{dijkgraaf} 

In the early 1990s (\cite{dw:tgtagc}), the Robbert Dijkgraaf noticed
that a TQFT approach yields 
a beautiful and elegant solution to a
classical mathematical problem:  counting ramified and unramified covers
of a topological surface.

Let $ (X,p_1,\ldots,p_r,q_1,\ldots,q_s) $ be an $(r+s)$-marked smooth topological surface.  Let $
\underline{\eta}=(\eta_1,\ldots,\eta_r) $ be a vector of partitions of the
integer $ d $.  We define the \textit{Hurwitz number}:
$$
    \begin{array}{ccc}
        H^{h,X}_d(\underline{\eta}) & := &\mbox{weighted number of}\ \left\{%
        \begin{array}{c}
            degree \ d \ covers \\
            C \stackrel{\pi}{\longrightarrow} X\ such\ that:  \\
            \bullet\ C\ is\ a\ surface\ of\ genus\  h;\ \ \ \ \ \ \ \ \ \ \ \\
            \bullet\ \pi \ is \ unramified\ over\ \ \ \ \ \ \ \ \ \ \ \ \ \ \ \ \ \ \\
             X\setminus\{p_1,\ldots, p_r,q_1,\ldots,q_s\};\\
            \bullet \ \pi\ ramifies\ with\ profile\ \eta_i\ over\ p_i; \\ 
            \bullet\ \pi\ has\ simple\ ramification\ over\ q_i. 
        \end{array}
        \right\}
    \end{array}
$$
The above number is weighted by the number of automorphisms of such
covers.

For a Hurwitz number to be nonzero, $s$, $h$ and $\underline{\eta}$ must satisfy the Riemann Hurwitz formula. This is why we omit $s$ from the above notation. In particular, if we require $s=0$, then (at most one value of) $h$ is determined by $\underline{\eta}$. We denote by $H^{X}_d(\underline{\eta})$ the corresponding Hurwitz number.  

\vspace{0.5cm}
We define the TQFT $ \mathcal{D} $ as follows:
\begin{enumerate}
    \item
        the ground field is  $ \mathbb{C} $;
    \item
        the Hilbert space is $ H=\bigoplus_{\eta \vdash d} \mathbb{C}\mbox
        {e}_\eta $;
    \item
        morphisms are assigned according to the prescription:
        $$
            \begin{array}{ccl}
               \parbox[c]{0.20\textwidth}{\includegraphics[height=3.5cm]{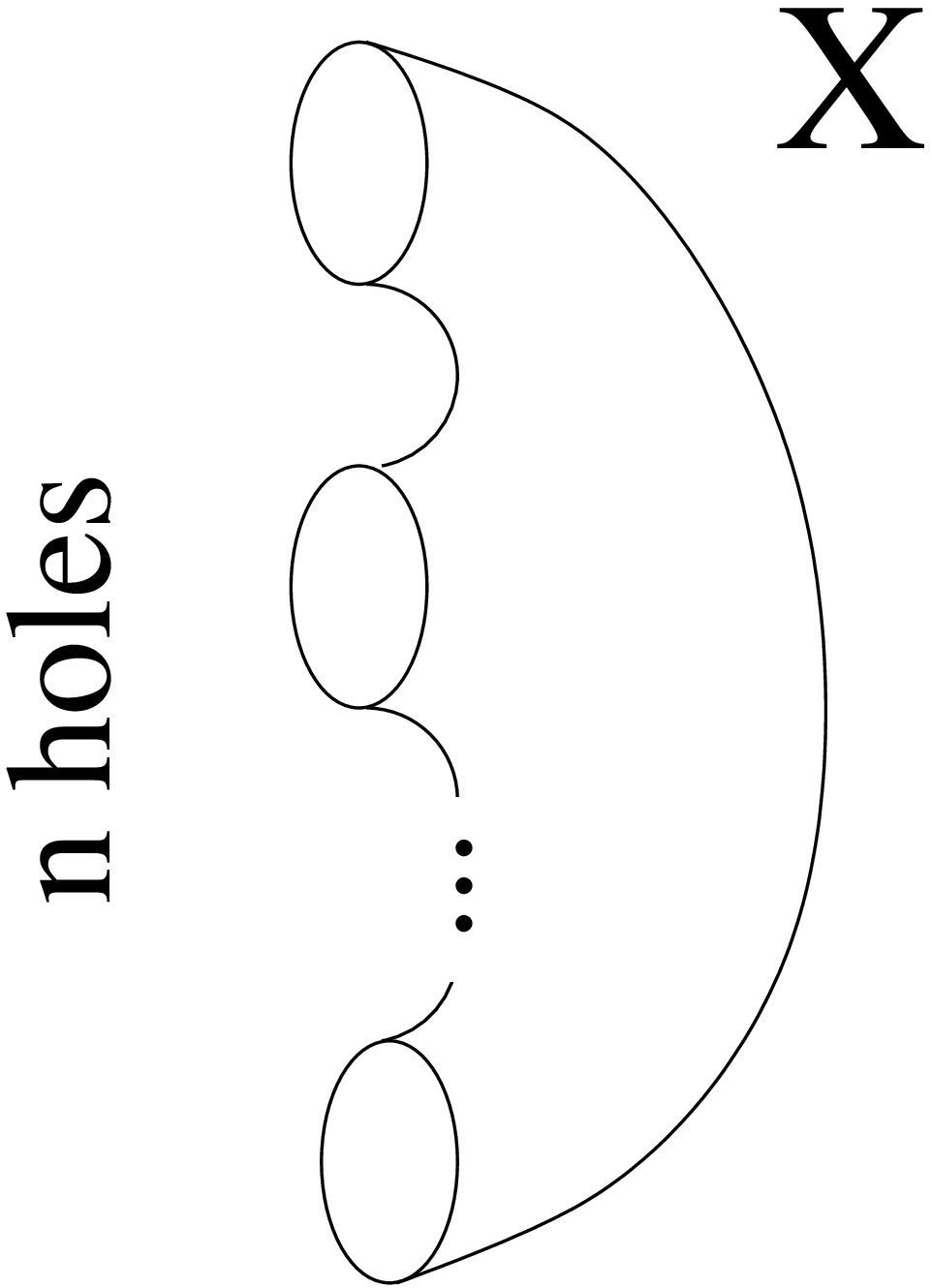}} & \stackrel{\mathcal
                {D}}{\mapsto} &
                \begin{array}{rcccl}
                    \mathcal{D}(X)& :& H^{\otimes n} & \longrightarrow &
                    \mathbb{C} \\
                    & & \mbox{e}_{\eta_1}\otimes \ldots\otimes \mbox{e}_{\eta_n}
                    & \mapsto & H^X_d(\underline{\eta}).
                \end{array}
                \\

              \parbox[c]{0.20\textwidth}{\includegraphics[height=2.7cm]{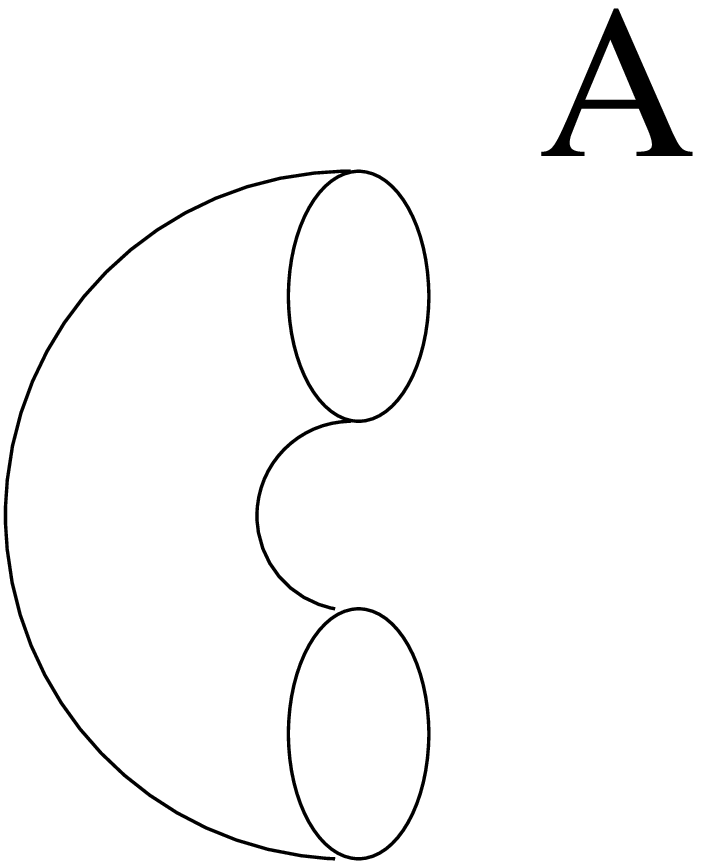}} & \stackrel{\mathcal
                {D}}{\mapsto} &\mathcal{D}(A)=\sum\mathfrak{z}(\eta)\mbox{e}_{\eta}\otimes
                \mbox{e}_{\eta}.
            \end{array}
        $$
\end{enumerate}

\begin{fact}[Dijkgraaf, Witten/ Freed, Quinn]
    The above assignment defines a semisimple TQFT $ \mathcal{D} $.
    Let $ \eta $ be a partition of $ d $, representing a conjugacy class
    of the symmetric group, and let $ h $ be an element in this
    conjugacy class. Via the identification:
    $$
        \mbox{e}_\eta= \frac{1}{d!}\sum_{g\in S_d}g^{-1}hg,
    $$
    the Hilbert space is isomorphic, as a Frobenius algebra, to the
    class algebra of the symmetric group in $ d $ letters, $ \mathcal{Z}
    (\mathbb{C}[S_d]) $.

    A semisimple basis is indexed by irreducible representations $ \rho $
    of $ S_d $.  Let $ \rho $ be such a representation and $ \mathcal{X}_\rho
    $ its character function, then:
    $$
        \mbox{e}_\rho=(dim\rho)\sum_{\eta\vdash d} \mathcal{X}_\rho(\eta)
        \mbox{e}_\eta.%
        \label{ssb}
    $$
\end{fact}

This allows Dijkgraaf to recover the classical Burnside formula,
expressing the number of unramified covers of a genus $ g $ curve as:
\begin{eqnarray}
    \begin{vimp}
        H_d^{gd-d+1,g}(\phi)=\displaystyle{\sum_{\rho} \left({\frac{d!}{\mbox{dim}\rho}}\right)^
        {2g-2}.}
    \end{vimp}
    \label{burn}
\end{eqnarray}

\subsection{Weighted TQFT's}
\label{emb} A weighted TQFT contains some extra structure with respect
to an ordinary TQFT\@.  Every cobordism comes equipped
with a sequence of weights, or levels.  When you concatenate two
cobordisms, you add the levels componentwise.  We are in particular
interested in the theory with $ 2 $ levels.

Define the category \textbf{2Cob$ ^{k_1,k_2} $} as follows:
\begin{enumerate}
    \item
        Objects and tensor structure are the same as in \textbf{2Cob}.
    \item
        Morphisms are given by triples $ (W,k_1,k_2) $, where $ W $ is
        an oriented cobordism as in \textbf{2Cob}, $ k_1,k_2 $ are two
        integers called levels.
    \item
        Composition of morphisms consists in concatenating the
        cobordisms and adding the levels componentwise.
\end{enumerate}
\begin{defi}
    A weighted TQFT is a functor of tensor categories:
    $$
        \mathcal{WT}:\mathbf{2Cob^{k_1,k_2}}\longrightarrow \mathbf{FRMod}.
    $$
\end{defi}

It is immediate  that if we restrict our attention to only
cobordisms with weight $ (0,0) $, we obtain an ordinary TQFT\@.  More
generally, there are a $ \mathbb{Z}\times\mathbb{Z} $ worth of ordinary
TQFTs embedded in a weighted TQFT\@.  Denote by $ \mathcal{X} $ the Euler
characteristic of a cobordism $ W $.  For any $ (a,b)\in\mathbb{Z}\times\mathbb
{Z} $, restricting the weighted TQFT to cobordisms with level
$$
    (a\mathcal{X},b\mathcal{X})
$$
yields an ordinary TQFT\@.

\subsection{Generation Results} 

\begin{figure}[htbp]
    \begin{center}
        \includegraphics[width=10cm]{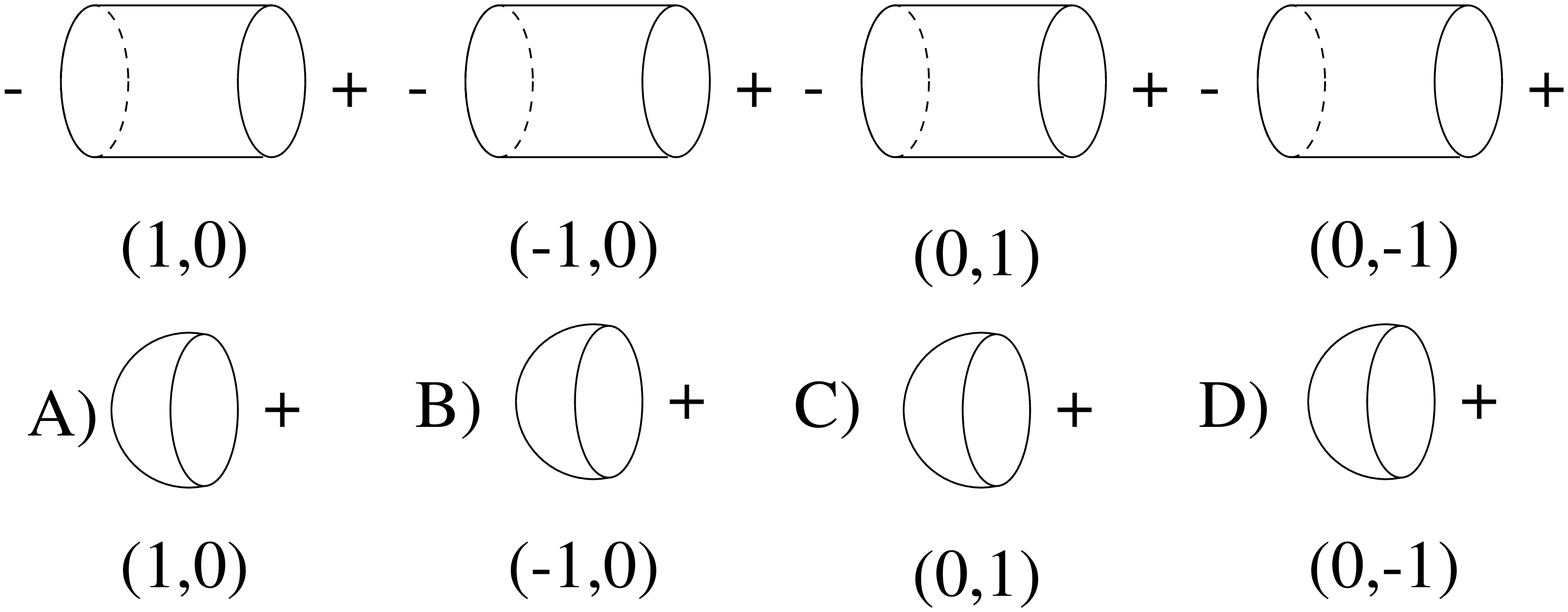}
        \caption{Level changing objects.}%
        \label{wt}
    \end{center}
\end{figure}

There are several possible ways to
generate a weighted TQFT\@.  A particularly natural one consists in
generating the level $ (0,0) $ TQFT, and then giving natural operators
that allow one to shift the levels. These elements can be chosen to be, for
example, the cylinders with weight $ (\pm1,0) $ and $ (0,\pm1) $.  These
operators change the levels of the cobordisms without altering its
topology.
An equivalent, and equally natural choice, is given by the caps, as illustrated in Figure~\ref{wt}.

In particular, it is immediate to see that A (resp. C) is the inverse of
B (resp. D) in the level $ (0,0) $ Frobenius algebra. Hence the following
generation result.
\begin{fact}[Bryan-Pandharipande, \cite{bp:tlgwtoc} 4.1] \label{wtqftgen}
    A weighted TQFT $ \mathcal{WT} $ is uniquely determined by a
    commutative Frobenius algebra over $ k $ for the level $ (0,0) $
    theory and by two distinguished invertible elements in the Frobenius
    algebra:
    $$
        \mathcal{WT}\left(
        \begin{array}{c}
            \includegraphics[height=1cm]{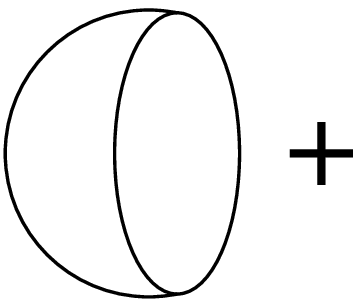}\\
            (-1,0)
        \end{array}
        \right),\ \ \mathcal{WT}\left(
        \begin{array}{c}
            \includegraphics[height=1cm]{unit.eps}\\
            (0,-1)
        \end{array}
        \right).
    $$
\end{fact}
\subsection{Semisimple Weighted TQFT}

\begin{figure}[htbp]
    $$
        \begin{array}{ccc}
            \begin{array}{c}
                \includegraphics[height=1cm]{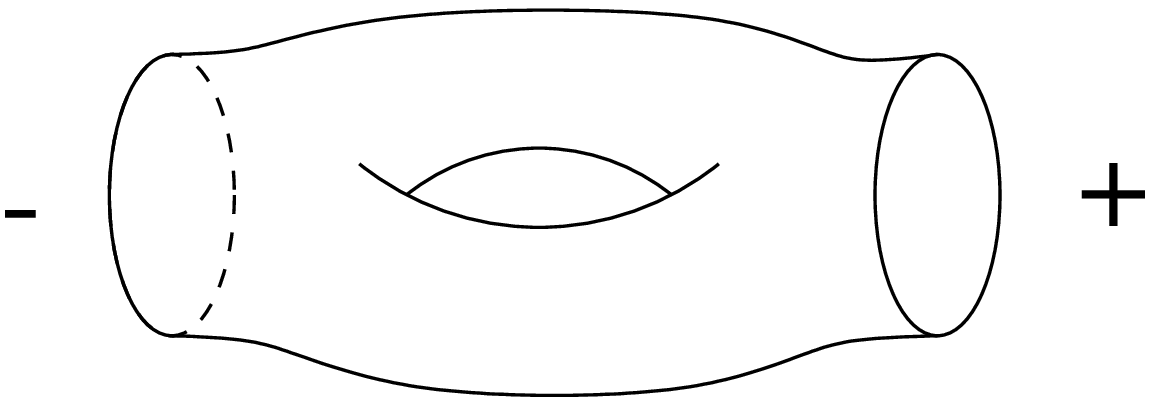}\\
                (0,0)
            \end{array}
            &
            \begin{array}{c}
                \includegraphics[height=1cm]{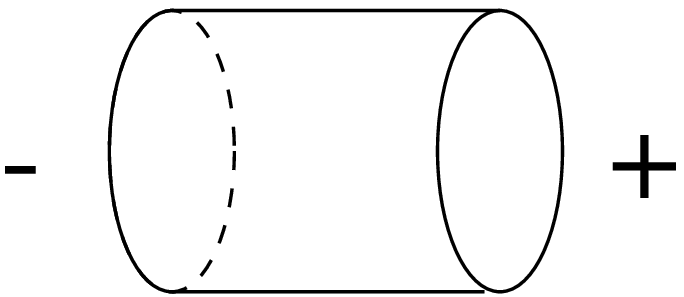}\\
                (-1,0)
            \end{array}
            &
            \begin{array}{c}
                \includegraphics[height=1cm]{cili.eps}\\
                (0,-1)
            \end{array}
            \\
            \downarrow & \downarrow & \downarrow \\
            \lambda_i & \mu_i & \overline{\mu}_i

        \end{array}
    $$
    \caption{The genus-adding and the level-changing operators.}%
    \label{ga}
\end{figure}

A weighted TQFT of rank $r$ is semisimple if there is a basis for the
Hilbert space such that all the non-zero tensors in the
theory are diagonal.  This is equivalent to asking that all embedded
ordinary TQFT's are semisimple (possibly with rescaled semisimple bases).
Let $ \lambda_1,\dots,\lambda_r $ be the eigenvalues of the level $ (0,0)
$ genus adding operator. Let $ \mu_1,\dots,\mu_r $ be the eigenvalues of
the level $ (-1,0) $ annulus, and $ \overline{\mu}_1,\dots,\overline{\mu}_r
$ be the eigenvalues for the level $ (0,-1) $ annulus, as illustrated in
Figure~\ref{ga}.
\begin{fact}[Bryan-Pandharpiande,\cite{bp:tlgwtoc}, 5.2]
    \label{wtqftstructure} Let $ \mathcal{WT} $ be a semisimple TQFT\@.
    Denote by $ W_m^n(g|k_1,k_2) $ a cobordism of genus $ g $ between $
    m $ input and $ n $ output holes, of level $ (k_1,k_2) $. Then:
    $$
        \mathcal{T}(W_m^n(g|k_1,k_2))= \sum_{i=1}^r\lambda_i^{g+n-1}\mu_i^
        {-k_1}\overline{\mu}_i^{-k_2}\underbrace{\mbox{e}^i\otimes\dots\otimes\mbox
        {e}^i}_{m\ %
        \ %
        times}\underbrace{\otimes\mbox{e}_i\otimes\dots\otimes\mbox{e}_i}_
        {n \ \ {times}}.
    $$
    In particular:
    \begin{eqnarray}
        \mathcal{T}(W_0^0(g|k_1,k_2))= \sum_{i=1}^r\lambda_i^{g-1}\mu_i^
        {-k_1}\overline{\mu}_i^{-k_2}.
    \end{eqnarray}
\end{fact}

\textbf{Observation:} the following equivalent definitions can be given
for the quantities $ \lambda_i $, $ \mu_i $ and $ \overline{\mu}_i $.
Denote by e$ _1,\dots $, e$ _r $ the vectors of a semisimple basis for
the weighted TQFT $ \mathcal{WT} $:
\begin{itemize}
    \item
        $ \lambda_i^{-1} $ is the value of the level $ (0,0) $ counit on
        e$ _i $:
        $$
            \mathcal{WT}\left(
            \begin{array}{c}
                \includegraphics[height=1cm]{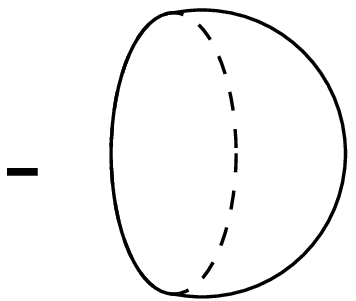}\\
                (0,0)
            \end{array}
            \right)(\mbox{e}_i)= \lambda_i^{-1}.
        $$
    \item
        $ \mu_i $ is the coefficient of e$ _i $ in the level $ (-1,0) $
        +disc vector:
        $$
            \mathcal{WT}\left(
            \begin{array}{c}
                \includegraphics[height=1cm]{unit.eps}\\
                (-1,0)
            \end{array}
            \right)= \sum \mu_i\mbox{e}_i.
        $$
    \item
        $ \overline{\mu}_i $ is the coefficient of e$ _i $ in the level $
        (0,-1) $ +disc vector:
        $$
            \mathcal{WT}\left(
            \begin{array}{c}
                \includegraphics[height=1cm]{unit.eps}\\
                (0,-1)
            \end{array}
            \right)= \sum \overline{\mu}_i\mbox{e}_i.
        $$
\end{itemize}

\section{Construction of the Theory}\label{u}
\subsection{The Admissible Covers Invariants}
\label{aci}

Let $(X,p_1,\ldots,p_r)$ be a smooth, irreducible,
 projective curve of genus $g$ with $r$ distinct marked points, and
$N=L_1\oplus L_2$ a rank $2$ vector bundle on $X$. The torus $T=\Cstar
\times \Cstar$ acts naturally on $N$:\label{toract} the first coordinate scales
(with weight one) the fibre of $L_1$, the second coordinate scales the
fibre of $L_2$.

The T-equivariant cohomology of a point is a polynomial ring in two indeterminates, that we denote:
$$H^{\ast}_T(pt)=\mathbb{C}[s_1,s_2].$$
We are interested in the following class of intersection numbers:
$$A_d^h(N):=\int_{ \overline{Adm}_{h\stackrel{d}{\rightarrow}X,
(\eta_1 p_1,\ldots,\eta_rp_r)}}
 \hspace{-1.5cm} e^{eq}(-R^\bullet\pi_\ast f^\ast(L_1\oplus L_2)),$$
where:
\begin{itemize}
  \item $ \overline{Adm}_{h\stackrel{d}{\rightarrow}X(\eta_1 p_1,\ldots,\eta_rp_r)}$
 is as defined on page \pageref{admcov}.
    \item $e^{eq}$ is the equivariant Euler class of the virtual bundle in question.
      \item $\pi$ is the universal family over the space of admissible covers.
        \item $f$ is the universal cover map followed by the canonical contraction map
         to $X$.
\end{itemize}

By \cite{b:bps}, this integral only depends on the genus $g$ of the curve $X$ and on
the degrees $k_1$ and $k_2$ of the line bundles $L_1$ and
$L_2$.  In our forthcoming
TQFT formulation it will be useful to emphasize this fact, so we
choose to denote the above invariants:
$$A_d^h(N)_{\underline{\eta}}=A_d^h(g|k_1,k_2)_{\underline{\eta}}.$$

We consider these invariants for all genera $h$, and organize them in generating function form as follows:
\begin{eqnarray}
A_d(g|k_1,k_2)_{\underline{\eta}}:= \sum_{h\in\mathbb{Z}}u^{\star(h)}A_d^h(g|k_1,k_2)_{\underline{\eta}}.
\label{adminv}
\end{eqnarray}

The appropriate exponent for the generating function is defined:
$$\star(h)= dim ( \overline{Adm}_{h\stackrel{d}{\rightarrow}X,
(\eta_1 p_1,\ldots,\eta_rp_r)}) = 2h-2 +d(2-2g-r) +\sum_{i=1}^r \ell(\eta_i) .$$

By expanding the equivariant Euler class in terms of ordinary Chern classes and of the equivariant parameters,
we can express these invariants in terms of nonequivariant integrals.
Let $h\in \mathbb{Z}\cup \phi$ be a function of $b_1, b_2$ determined by the equation
$$b_1+b_2 =dim ( \overline{Adm}_{h\stackrel{d}{\rightarrow}X,(\eta_1x_1,\cdots,\eta_rx_r)}) = 2h-2 +d(2-2g-r)+ \sum_{i=1}^r \ell(\eta_i) .$$
Define
$$A_d^{b_1,b_2}(g|k_1,k_2)_{\underline{\eta}}:=
\int_{\overline{Adm}_{h\stackrel{d}{\rightarrow}X,(\eta_1x_1,\cdots,\eta_rx_r)}}\hspace{-1.5cm}c_{b_1}(-R^\bullet\pi_\ast
f^\ast(L_1))c_{b_2}(-R^\bullet\pi_\ast f^\ast(L_2)).$$
Then the relative invariants are:
\begin{eqnarray}
A_d(g|k_1,k_2)_{\underline{\eta}}:= \sum_{b_1+b_2=0}^\infty u^{b_1+b_2 } s_1^{r_1-b_1}s_2^{r_2-b_2}A_d^{b_1,b_2}(g|k_1,k_2)_{\underline{\eta}}.\label{ordchern}\end{eqnarray}

This shows that the partition function for our invariants is a Taylor
series in $u$, whose coefficients are rational functions in $s_1$ and
$s_2$. The degree of these rational functions is
independent of $h$. It is equal to
$r_1+r_2-b_1-b_2=d(2g-2-r) -\sum_{i=1}^r \ell(\eta_i) .$

\subsection{The Weighted TQFT $\mathcal{U}$}
We construct a weighted TQFT\ \  $\mathcal{U}$, whose structure coefficients encode the invariants just presented. 

The ground ring is defined to be $R= \mathbb{C}[[u]](s_1,s_2)$.

The Hilbert space of the theory is a free $R$-module of rank equal to the number of partitions of the integer $d$. A privileged basis is  indexed by such partitions $\eta$.
$$\mathcal{H}= \bigoplus_{{\eta} \vdash d} R\mbox{e}_\eta.$$
We  denote the dual space by $\mathcal{H^\ast}$, and the dual basis vectors by e$^\eta$.

In order to  construct our TQFT we reason topologically. We think of the marked points on a curve $(X, x_1,\dots,x_{r+s})$ as of punctures that we can ``enlarge'' into  loops. We can assign positive or negative orientation to such loops, and arrange the negatively oriented loops  $x_1,\dots,x_r$ to the left, the positively oriented to the right (after relabelling $x_{r+i}=y_i$). We now have an oriented cobordism.

To completely determine the structure of the theory we  define the scalar maps associated to arbitrary cobordisms into the empty set, and the coproduct, that allows us to ``move'' boundary components from the left to the right.
$$\begin{array}{ccl}
\parbox[c]{0.35\textwidth}{\includegraphics[height=3.5cm]{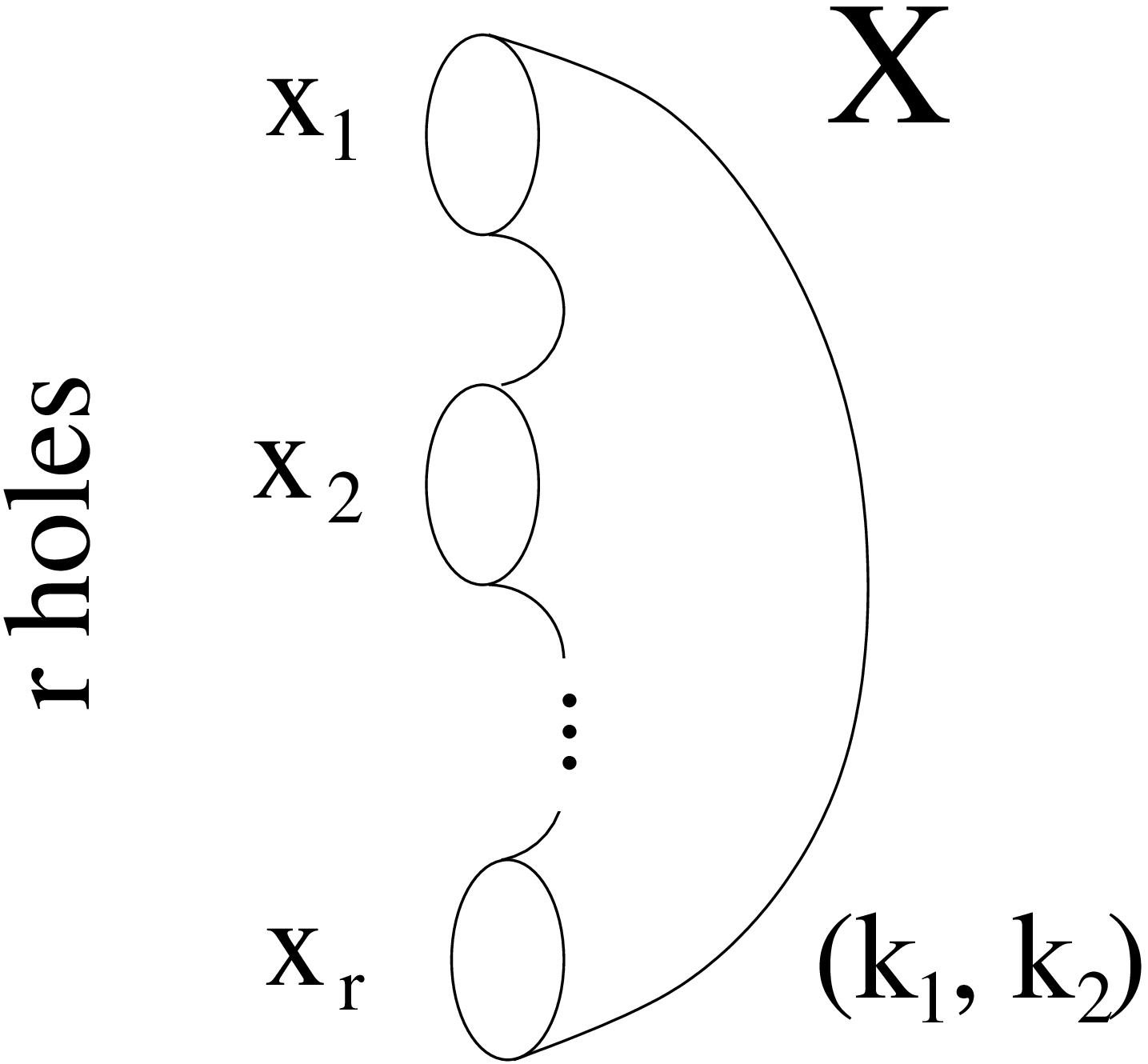}}
  & \stackrel{\mathcal{U}}{\mapsto} & \begin{array}{rcccl}
\mathcal{U}(X)& :& H^{\otimes r} &  \longrightarrow & \mathbb{C}[[u]](s_1,s_2) \\
              &  &               &                  &                          \\
              &  & \mbox{e}_{\eta_1}\otimes ...\otimes \mbox{e}_{\eta_r} &
\mapsto &  A_d(g|k_1,k_2)_{\underline{\eta}} .
\end{array}\\
 & & \\

\parbox[c]{0.35\textwidth}{\includegraphics[height=2.7cm]{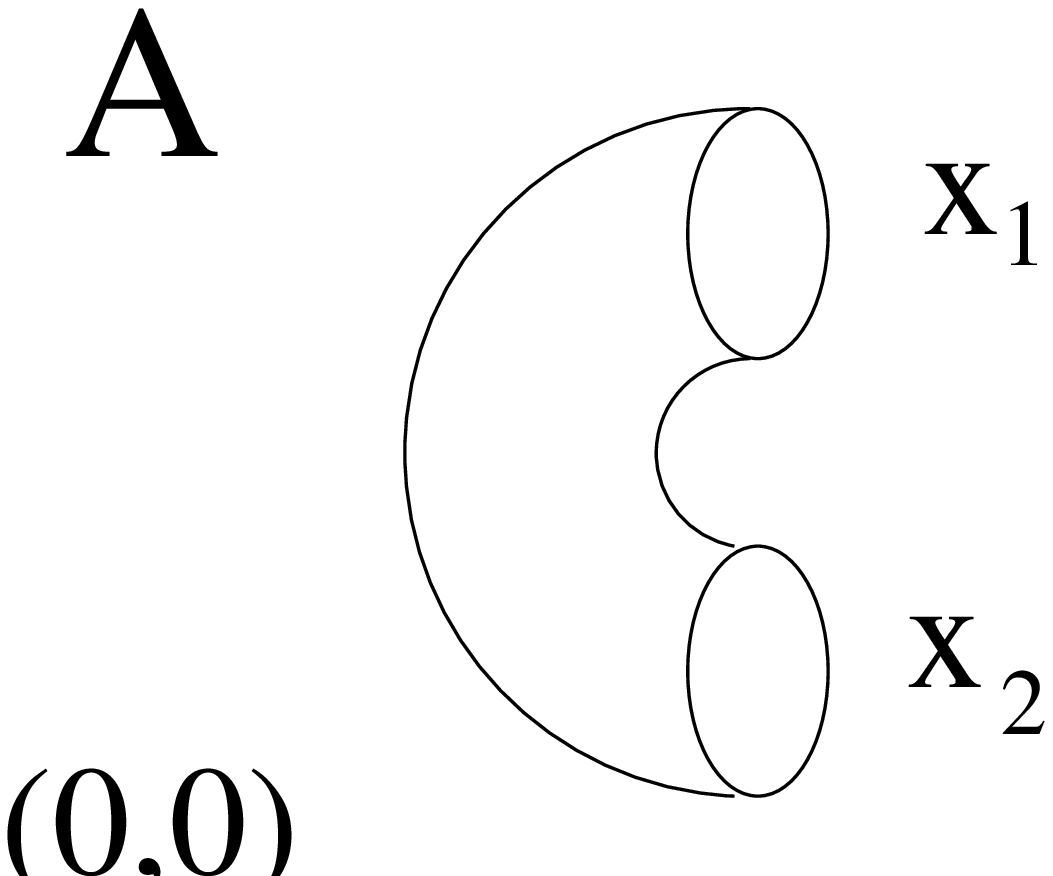}}
 & \stackrel{\mathcal{U}}{\mapsto}
&\mathcal{U}(A)=\sum_{{\eta} \vdash d}\mathfrak{z}(\eta)(s_1s_2)^{\ell(\eta)}\mbox{e}_{\eta}\otimes
\mbox{e}_{\eta}.
\end{array}
$$

The combinatorial factor $\mathfrak{z}(\eta)$ is defined in page \pageref{zeta}.
\begin{theorem}
The structure $\mathcal{U}$ defined in the previous paragraph is a two-level, weighted semisimple TQFT.
\label{tqftteor}
\end{theorem}

In practical terms, it is often very convenient to
adopt the conventional riemannian geometry tensor notation. If
$X=(X,x_1,\dots,x_r,y_1,\dots,y_s|k_1,k_2)$ represents a cobordism of genus
$g$ and level $k_1,k_2$ from $r$ circles to $s$ circles, then $\mathcal{U}(X)$ is an element of $({\mathcal{H}^\ast})^{\otimes r}\otimes {\mathcal{H}}^{\otimes s}$.  We denote by

$$ A_d(g|k_1,k_2)_{{\eta_1,\dots,\eta_r}}^{{\mu_1,\dots,\mu_s}}$$
the coordinate of $\mathcal{U}(X)$ in the direction of the basis element  e$^{\eta_1}\otimes \dots \otimes$ e$^{\eta_r}\otimes$e$_{\mu_1}\otimes \dots \otimes$ e$_{\mu_s}$.

With this notation, the coproduct gives the following formula for raising and lowering indices:
\begin{eqnarray}A_d(g|k_1,k_2)_{{\eta_1,\dots,\eta_r}}^{{\mu_1,\dots,\mu_s}}= \left( \prod_{i=1}^s \mathfrak{z}(\mu_i) (s_1s_2)^{\ell(\mu_i)} \right)
 A_d(g|k_1,k_2)_{{\eta_1,\dots,\eta_r,\mu_1,\dots,\mu_s}}.\label{indices}\end{eqnarray}

\subsection{Proof of Theorem \ref{tqftteor}}

Proving that $\mathcal{U}$ is indeed a  weighted TQFT 
amounts to verifying the following three statements:
\begin{description}
  \item[identity:] the tensor associated to the level $(0,0)$ trivial cobordism from the circle to the circle is the identity morphism of the Hilbert space $\mathcal{H}$.
    \item[gluing two curves:] for any two vectors $\underline{\eta},\underline{\mu}$ of partitions of $d$, and integers satisfying $g=g'+g''$, $ k_1=k_1'+k_1''$, $ k_2=k_2'+k_2''$,
\begin{eqnarray}A_d(g|k_1,k_2)_{{\eta_1,\dots,\eta_r}}^{{\mu_1,\dots,\mu_s}}= 
\sum_{\nu \vdash d}A_d(g'|k'_1,k'_2)_{{\eta_1,\dots,\eta_r}}^{\nu}
 A_d(g''|k''_1,k''_2)_{\nu}^{{\mu_1,\dots,\mu_s}}. \label{gtc}\end{eqnarray}
      \item[self-gluing:] for any vector of partitions $\underline{\eta}$, and integers $g$, $k_1$, $k_2$,
\begin{eqnarray}A_d(g+1|k_1,k_2)_{{\eta_1,\dots,\eta_r}} =
  \sum_{\nu \vdash d}A_d(g|k_1,k_2)_{{\eta_1,\dots,\eta_r,\nu}}^\nu.\label{selfg}\end{eqnarray}
\end{description}

\subsubsection{Identity}

This fact is easily proven. One very clever way to do it,
 which is pursued in \cite{bp:tqft}, is to notice that the degree $0$
 coefficients in our TQFT agree with the classical TQFT of Hurwitz numbers
 constructed by Dijkgraaf in \cite{dw:tgtagc} and recalled in section \ref{dijkgraaf}.
 The vanishing of all higher degree terms can be obtained as a straightforward
 consequence of the gluing laws, or simply by showing that the dimensions of the
moduli spaces in question exceed the maximum degree of a non-equivariant class in the integrand.

\subsubsection{Gluing Two Curves}

In order to minimize the burden of bookkeeping, we  prove the result when  $r=s=0$  (i.e.\,, the resulting glued curve is not marked). In the general case, the proof follows exactly the same steps, and all the extra indices  are simply carried along for the ride.

Consider a one parameter family of genus $g$ curves $W$, and the corresponding  map to the moduli space,
$$
\begin{array}{cccc}
 & W & & \\
 &\downarrow & & \\
\varphi: & \mathbb{A}^1 & \rightarrow  & \overline{M}_g,
\end{array}$$
such that:
\begin{itemize}
  \item the central fiber
$$W_0 = X_1 \displaystyle{\bigcup_{b_1=b_2}} X_2$$
is a nodal curve obtained by attaching at a point two smooth curves of genus $g'$ and $g''$ (with $g' +g''=g$);
\item all other fibers $W_s, s\not=0$, are smooth curves of genus $g$.
\end{itemize}
Consider the moduli space $\overline{Adm}^{}_{h\stackrel{d}{\rightarrow}g}$ of admissible covers of a genus $g$ curve by a genus $h$ curve, all ramification simple. By \cite{acv:ac}, there is a  flat morphism
$$\overline{Adm}^{}_{h\stackrel{d}{\rightarrow}g}\rightarrow\overline{M}_{g},$$

We can construct the following cartesian diagram:
$$
\begin{array}{ccccc}
\mathcal{A}_s= \overline{Adm}^{}_{h\stackrel{d}{\rightarrow}W_s} & \hookrightarrow & \mathcal{A} & \rightarrow & \overline{Adm}^{}_{h\stackrel{d}{\rightarrow}g} \\
 & & & & \\
\downarrow & &\downarrow & & \downarrow \\
 & & & & \\
\{s\} & \hookrightarrow & \mathbb{A}^1 & \rightarrow & \overline{M}_g
\end{array}
$$

The stack $\mathcal{A}$ must be thought as of the stack of relative admissible covers of the family $W$. For $s\not=0$, we obtain admissible covers of a smooth genus $g$ curve; for $s=0$, we recover admissible covers of the nodal curve $W_0$.

It is  possible to construct two line bundles $\mathcal{L}_1$ and $\mathcal{L}_2$ on $W$, with the following properties:
\begin{enumerate}
  \item $\mathcal{L}_i$ restricted to any  fiber $W_s$ is a line bundle $L_{i,s}$ of degree $k_i$.
    \item Over the central fiber $W_0$,  $\mathcal{L}_i$ restricts to a line bundle $L'_{i,s}$ of degree $k_i'$ on $X_1$, and restricts to a line bundle $L''_{i,s}$ of degree $k_i''$ on $X_2$.
      \item $\Cstar$ acts naturally on $\mathcal{L}_i$ by scaling the fibers (with weight one).
\end{enumerate}
Consider the following diagram:
$$
\begin{array}{ccccc}
\mathcal{U_A} & \stackrel{f}{\rightarrow} & \mathcal{W} & \rightarrow & W \\
 & & & & \\
\pi \downarrow &  \swarrow & & &  \\
 & & & &  \\
\mathcal{A} & &  & &
\end{array}
$$

$\mathcal{U_A}$ is the universal family of the moduli space $\mathcal{A}$,  $\mathcal{W}$ is the universal target and $f$ the universal admissible cover map.

The pull-push
$$\mathcal{I} = -R^\bullet\pi_\ast f^\ast (\mathcal{L}_1\oplus\mathcal{L}_2)$$
is a virtual bundle of virtual rank $r=2g-2-d(k_1+k_2)$.

By the flatness of the family $\mathcal{A}$ over $\mathbb{A}^1$, the integral of the top Chern class $c_r(\mathcal{I})$ restricted to a fiber $\mathcal{A}_s$ is independent of the fiber. For $s\not=0$, we obtain
\begin{eqnarray}
\int_{\overline{Adm}^{}_{h\stackrel{d}{\rightarrow}W_s}}\hspace{-0.5cm}c_r(\mathcal{I}\mid_s)= A^h_d(g|k_1,k_2)
\end{eqnarray}

We  want to evaluate the same expression restricted to s=0, and show it equals the right hand side
of (\ref{gtc}). We choose to show the equality at the generic genus $h$ degree of the generating function,
to emphasize the geometric nature of the construction.
We hence need to establish the following claim, which consists of expanding the genus $h$ term in equation (\ref{gtc}), and lowering indices as in (\ref{indices}).
\begin{Claim}
$$\int_{\overline{Adm}^{}_{h\stackrel{d}{\rightarrow}W_0}}\hspace{-0.5cm}c_r(\mathcal{I}\mid_0)=
\sum_{\nu\vdash d} \mathfrak{z}(\nu) (s_1s_2)^{\ell(\nu)}\sum_{h_1,h_2}A^{h_1}_d(g'|k'_1,k'_2)_\nu A^{h_2}_d(g''|k''_1,k''_2)_\nu,$$
where the second sum is over pairs of indices such that  $h_1 +h_2 +\ell(\nu)-1= h $.
\end{Claim}

\textsc{Proof:} Recall that, by  (\ref{rnodal})  in section~\ref{nodal},
$$[\overline{Adm}^{}_{h\stackrel{d}{\rightarrow}W_0}]= \sum_{\nu \vdash d} \mathfrak{z}(\nu)\sum_{h_1,h_2}[\overline{Adm}^{}_{h_1\stackrel{d}{\rightarrow}X_1,(\nu b_1)}]\times[\overline{Adm}^{}_{h_2\stackrel{d}{\rightarrow}X_2,(\nu b_2)}],$$
where:
\begin{itemize}
  \item $h_1 +h_2 +\ell(\nu)-1= h$;
    \item $dim(\overline{Adm}^{}_{h_1\stackrel{d}{\rightarrow}X_1,(\nu b_1)}) + dim(\overline{Adm}^{}_{h_2\stackrel{d}{\rightarrow}X_2,(\nu b_2)}) = dim(\overline{Adm}^{}_{h\stackrel{d}{\rightarrow}W_0})  $.
\end{itemize}

Consider the pull-back of the normalization sequence associated to the restriction of $\mathcal{L}_i$ to $W_0$:
$$0\rightarrow f^\ast(L_{i,0})\rightarrow f^\ast(L'_{i,0})\oplus f^\ast(L''_{i,0}) \rightarrow f^\ast(L_{i,0})\mid_{X_1\cap X_2}\rightarrow 0.$$

This sequence yields a long exact sequence of higher direct image sheaves
$$0\rightarrow R^0\pi_\ast f^\ast (L_{i,0})\rightarrow R^0\pi_\ast f^\ast (L'_{i,0})\oplus R^0\pi_\ast f^\ast(L''_{i,0})\rightarrow R^0\pi_\ast f^\ast(L_{i,0})\mid_{X_1\cap X_2}\rightarrow $$
$$\rightarrow R^1\pi_\ast f^\ast (L_{i,0})\rightarrow R^1\pi_\ast f^\ast (L'_{i,0})\oplus R^1\pi_\ast f^\ast(L''_{i,0})\rightarrow 0.$$

Notice that $(L_{i,0})\mid_{X_1\cap X_2}$ is  a skyscraper sheaf $\mathbb{C}_b$, on which $\Cstar$ acts with weight $1$. 

Let us  restrict our attention to a connected component of $\mathcal{A}_0$ on which the covers split as two smooth covers of genus $h_1$ and $h_2$, with ramification profile $\nu$ over the shadows of the node. Here, $f^\ast (L_{i,0})\mid_{X_1\cap X_2}$ is a trivial vector bundle of rank $\ell(\nu)$, endowed with a natural $\Cstar$ action.

From this fact and the exact sequence above we conclude:
$$c_{r_i}( -R^\bullet \pi_\ast f^\ast (L_{i,0}))= s_i^{\ell(\nu)}c_{r'_i}( -R^\bullet \pi_\ast f^\ast (L'_{i,0}))c_{r''_i}( -R^\bullet \pi_\ast f^\ast (L''_{i,0})),$$
and finally
\begin{eqnarray}
c_{r} (\mathcal{I}\mid_0)= (s_1s_2)^{\ell(\nu)}c_{r'}(\mathcal{I}\mid_0')c_{r''} (\mathcal{I}\mid_0'').
\end{eqnarray}

Now clinching the claim is just a matter of putting everything together:
$$\
\begin{array}{lcl}\hspace{-1cm}\displaystyle{
\int_{\overline{Adm}^{}_{h\stackrel{d}{\rightarrow}W_0}}\hspace{-0.5cm}c_r(\mathcal{I}\mid_0)} & = & \displaystyle{
 \sum_\nu \mathfrak{z}(\nu) \sum_{h_1,h_2}\int_{\overline{Adm}^{}_{h_1\stackrel{d}{\rightarrow}X_1,(\nu b_1)}\times\overline{Adm}^{}_{h_2\stackrel{d}{\rightarrow}X_2,(\nu b_2)}}\hspace{-0.5cm}c_{r}(\mathcal{I}\mid_0)} \\ & & \\ & & \\
  & = & \displaystyle{
\sum_\nu \mathfrak{z}(\nu)(s_1s_2)^{\ell(\nu)} \sum_{h_1,h_2}\int_{\overline{Adm}^{}_{h_1\stackrel{d}{\rightarrow}X_1,(\nu b_1)}}\hspace{-0.5cm}c_{r'} (\mathcal{I}\mid_0')\int_{\overline{Adm}^{}_{h_2\stackrel{d}{\rightarrow}X_2,(\nu b_2)}}\hspace{-0.5cm}c_{r''} (\mathcal{I}\mid_0'' )} \\  & & \\ & & \\

 &=& \displaystyle{
\sum_\nu \mathfrak{z}(\nu) (s_1s_2)^{\ell(\nu)}\sum_{h_1,h_2}A^{h_1}_d(g'|k'_1,k'_2)_\nu A^{h_2}_d(g''|k''_1,k''_2)_\nu}\ \ .

\end{array}
$$

\subsubsection{Self-gluing}
The structure of the proof is very similar to the previous case.
 Again, we simplify the notation by assuming $r=0$.

Consider a one parameter family of genus $g$ curves $W$, and the corresponding  map into the moduli space,
$$
\begin{array}{cccc}
 & W & & \\
 &\downarrow & & \\
\varphi: & \mathbb{A}^1 & \rightarrow  & \overline{M}_g,
\end{array}$$
such that:
\begin{itemize}
  \item the central fiber
$$W_0 = X/\{b_1=b_2\}$$
is a nodal curve obtained by identifying two distinct points on an irreducible smooth curve $X$ of genus $g-1$;
\item all other fibers $W_s, s\not=0$, are smooth curves of genus $g$.
\end{itemize}

As before, we can construct the following cartesian diagrams:
$$
\begin{array}{ccccc}
\mathcal{A}_s= \overline{Adm}^{}_{h\stackrel{d}{\rightarrow}W_s} & \hookrightarrow & \mathcal{A} & \rightarrow & \overline{Adm}^{}_{h\stackrel{d}{\rightarrow}g} \\
 & & & & \\
\downarrow & &\downarrow & & \downarrow \\
 & & & & \\
\{s\} & \hookrightarrow & \mathbb{A}^1 & \rightarrow & \overline{M}_g
\end{array}
$$
and  two line bundles $\mathcal{L}_1$ and $\mathcal{L}_2$ on $W$, with the following properties:
\begin{enumerate}
  \item $\mathcal{L}_i$ restricted to any  fiber $W_s$ is a line bundle $L_{i,s}$ of degree $k_i$.
    \item Over the central fiber $W_0$,  $\mathcal{L}_i$ pulls back to a line bundle $ L'_{i,s}$ of degree $k_i$ on the normalization $X$.
      \item $\Cstar$ acts naturally on $\mathcal{L}_i$ by scaling the fibers (with weight one).
\end{enumerate}

We now consider the equivariant top Chern class of the  pull-push
$$\mathcal{I} = -R^\bullet\pi_\ast f^\ast (\mathcal{L}_1\oplus\mathcal{L}_2).$$

For $s\not= 0$,
\begin{eqnarray}
\int_{\overline{Adm}^{}_{h\stackrel{d}{\rightarrow}W_s}}\hspace{-0.5cm}c_r(\mathcal{I}\mid_s)= A^h_d(g|k_1,k_2).
\end{eqnarray}

Again, we can show that the corresponding integral over the central fiber yields exactly the genus $h$ expansion of the right hand side
of equation (\ref{selfg}).

\begin{Claim}
$$\int_{\overline{Adm}^{}_{h\stackrel{d}{\rightarrow}W_0,}}\hspace{-0.5cm}c_r(\mathcal{I}\mid_0)=
\sum_\nu \mathfrak{z}(\nu) (s_1s_2)^{\ell(\nu)}A^{h'}_d(g-1|k_1,k_2)_{\nu,\nu} \ \ ,$$
where   $h'+\ell(\nu)= h $.
\end{Claim}

\textsc{Proof:} By  (\ref{irrnodal}) in section~\ref{nodal},
$$[\overline{Adm}^{}_{h\stackrel{d}{\rightarrow}W_0}]= \sum_{\nu \vdash d} \mathfrak{z}(\nu)[\overline{Adm}^{}_{h'\stackrel{d}{\rightarrow}X,(\nu b_1, \nu b_2)}],$$
with $h'+\ell(\nu)=h$

 As in the previous paragraph, after chasing the normalization sequence for the curve $W_0$, we obtain, over a connected component of $\mathcal{A}_0$ characterized by covers with ramification profile $\nu$ over the shadows of the node, the following decomposition:
\begin{eqnarray}
c_{r} (\mathcal{I}\mid_0)= (s_1s_2)^{\ell(\nu)}c_{r'}(\mathcal{I}\mid_0').
\end{eqnarray}

With this in hand, it is easy to conclude:
$$
\begin{array}{lcl}
\displaystyle{\int_{\overline{Adm}^{}_{h\stackrel{d}{\rightarrow}W_0}}\hspace{-0.5cm}c_r(\mathcal{I}\mid_0)} & = &
\displaystyle{\sum_\nu \mathfrak{z}(\nu) \int_{\overline{Adm}^{}_{h'\stackrel{d}{\rightarrow}X,(\nu b_1,\nu b_2)}}\hspace{-0.5cm}c_{r}(\mathcal{I}\mid_0)} \\ & & \\ & & \\
 & = & \displaystyle{\sum_\nu \mathfrak{z}(\nu)(s_1s_2)^{\ell(\nu)} \int_{\overline{Adm}^{}_{h'\stackrel{d}{\rightarrow}X,(\nu b_1,\nu b_2)}}\hspace{-0.5cm}c_{r'} (\mathcal{I}\mid_0')} \\ & & \\ & & \\

 & = & \displaystyle{\sum_\nu \mathfrak{z}(\nu) (s_1s_2)^{\ell(\nu)}A^{h'}_d(g-1|k'_1,k'_2)_{\nu,\nu}}\ \ .
\end{array}  $$

\section{Computing the Theory}\label{compute}

In order to determine the whole weighted TQFT it is sufficient to compute a small number of invariants, as seen in Fact~\ref{wtqftgen}. There are many possible choices for a set of generators; we make the following choice:\\
\textbf{generators for the level $(0,0)$ TQFT:}

\begin{enumerate}
    \item the coefficients $A_d(0|0,0)_\eta$ of the open (-)disc.
      \item the coefficients $A_d(0|0,0)^{\eta,\mu}$ of the (+,+) annulus.
      \item the coefficients $A_d(0|0,0)_{\eta,\mu,\nu}$ associated to the (-,-,-) pair of pants.
        \end{enumerate}
\textbf{generators for level shifting :}\\

4. the coefficients of the Calabi-Yau caps $A_d(0|-1,0)_\eta$ and $A_d(0|0,-1)_\eta$.
\vspace{0.5cm}

\begin{theorem}\label{equals}The level $(0,0)$ TQFT coincides with the level $(0,0)$ theory of Bryan and Pandharipande in \cite{bp:tlgwtoc}.
\end{theorem}

\textsc{Proof:} It is simple to compute independently the coefficients for 
the cap. Dimension counts show they are  degenerate, in the sense that only the constant term of 
the series is nonzero. The coefficients for the $(+,+)$-cylinder agree by definition. 
In section~\ref{pair} we conclude the proof by showing that the coefficients for the pair of pants are the
 same.

\vspace{0.5cm}
The significant difference in the theories lies in the Calabi-Yau caps.
These are computed by localization on the moduli spaces of admissible covers
in section \ref{CYcapsec}.

\subsection{The Level $(0,0)$ Pair of Pants}
\label{pair}

The invariants  $A^{\circ}_d(0|0,0)_{\eta,\nu,\mu}$ of the level $(0,0)$ pair  of pants are computed by the integrals:
$$\int_{\overline{Adm} ^\circ_{h\stackrel{d}{\rightarrow}\proj,(\eta 0, \mu 1, \nu\infty)}}\hspace{-1cm} c_{2h-2}^{eq} (-R^\bullet\pi_\ast f^\ast(\mathcal{O}_{\proj}\oplus \mathcal{O}_{\proj})).$$
The dimension of the moduli space in question is
$$2h -d -2 + \ell(\eta) +\ell(\mu) +\ell(\nu).$$
Hence, if $\ell(\eta) +\ell(\mu) +\ell(\nu) > d+2$, the relative connected integrals vanish. The disconnected integrals are then obtained inductively from  invariants of lower degree $d$.

All other invariants have contributions from connected components, and hence need to be computed directly.

In an appendix to \cite{bp:tlgwtoc}, Bryan, Pandharipande and Faber show that all invariants can be recursively determined from $A_d(0|0,0)_{(d),(d),(2)}$, the invariant corresponding to full ramification over two points, and a simple transposition over the third point. Their proof uses only TQFT formalism; hence it suffices to prove the following statement.

\begin{lemma}\label{pair}
For d$\geq$ 2,
$$A_d(0|0,0)_{(d),(d),(2)}= \frac{1}{2}\frac{s_1+s_2}{s_1s_2}\left(d\cot\left(\frac{du}{2}\right)-\cot\left(\frac{ u}{2}\right)\right).$$
\end{lemma}

\textbf{Note:} The above result differs from the analogous one in \cite{bp:tlgwtoc} by a factor of $-i$, that reflects a different normalization in their generating function conventions, that we do not wish to adopt.

\textsc{Proof:} Notice, first of all, that the full ramification conditions force our covers to be connected. In this case the connected and disconnected invariants coincide.

According to (\ref{ordchern}), we have:
$$A_d(0|0,0)_{(d),(d),(2)}= \sum_{b_1+b_2=0}^\infty u^{b_1+b_2 } s_1^{h-1-b_1}s_2^{h-1-b_2}\int_{\overline{Adm}_{h\stackrel{d}{\rightarrow}\proj,((d) 0, (d) 1, (2)\infty)}}\hspace{-2cm}c_{b_1}(\mathbb{E}^\ast)c_{b_2}(\mathbb{E}^\ast),$$
with $b_1+b_2$ equal to the dimension of the moduli space, which is
$$dim(\overline{Adm}_{h\stackrel{d}{\rightarrow}\proj,((d) 0, (d) 1, (2)\infty)})=2h-1.$$
For a given value of $h$, the only nonvanishing terms in the above expression are given by:
\begin{itemize}
  \item $b_1=h,\ \ \  b_2=h-1$;
  \item $b_1=h-1,\ \ \  b_2=h$.
\end{itemize}
Adding the two terms, we obtain
$$A_d^h(0|0,0)_{(d),(d),(2)}=\frac{s_1+s_2}{s_1s_2} \int_{\overline{Adm}_{h\stackrel{d}{\rightarrow}\proj,((d) 0, (d) 1, (2)\infty)}} \hspace{-2cm} -\lambda_h\lambda_{h-1}$$
and consequently, the generating function:
$$A_d(0|0,0)_{(d),(d),(2)}= \frac{s_1+s_2}{s_1s_2}\sum_{h=0}^\infty u^{2h-1 } \int_{\overline{Adm}_{h\stackrel{d}{\rightarrow}\proj,((d) 0, (d) 1, (2)\infty)}} \hspace{-2cm} -\lambda_h\lambda_{h-1} \ \ \ \ \ \   ,$$
where $\lambda_k$ denotes the $k^{th}$ Chern class of the (pull-back of the) Hodge bundle $\mathbb{E}$.

Let us recall that we defined the $\lambda$ classes on moduli spaces of admissible covers simply by pulling them back from the appropriate moduli spaces of stable curves. In particular observe the diagram:
$$
\begin{array}{ccc}

\overline{Adm}_{h\stackrel{d}{\rightarrow}\proj,((d) 0, (d) 1, (2)\infty)} & \stackrel{\rho}{\longrightarrow} & \overline{M}_{h,2} \\ & & \\
 &\searrow  & \downarrow \pi \\ & & \\
         & & \overline{M}_h
\end{array}
$$

The map $\rho$ is defined by marking on the admissible covers the unique preimages of the branch points $0$ and $1$. The Hodge bundle on $\overline{M}_h$ pulls back to the Hodge bundle on $\overline{M}_{h,2}$, hence we can think of the $\lambda$ classes on the moduli space of admissible covers as pulled back from $\overline{M}_{h,2}$.

Denote by $H_d\subset M_{h,2}$ the locus of curves admitting a degree $d$ map to $\proj$ which is totally ramified at the marked points. Let
$$\overline{H}_d\subset \overline{M}_{h,2}$$
be the closure of $H_d$, consisting of possibly nodal curves admitting a degree $d$ map to a tree of rational curves, fully ramified over the two marked points.
The image of the  map
$$\rho:\overline{Adm}_{h\stackrel{d}{\rightarrow}\proj,((d) 0, (d) 1, (2)\infty)} \longrightarrow \overline{M}_{h,2}$$
is precisely $\overline{H}_d$, and $\rho$ is a degree $2h$ map onto its image.

From this we conclude that
$$\int_{\overline{Adm}_{h\stackrel{d}{\rightarrow}\proj,((d) 0, (d) 1, (2)\infty)}} \hspace{-2cm} -\lambda_h\lambda_{h-1}\ \ \ \ \ \ =2h\int_{[\overline{H}_d]}  -\lambda_h\lambda_{h-1}.$$

This is exactly the  integral computed in \cite{bp:tlgwtoc}, pages 28-29, hence the result follows. This concludes the proof of  Lemma  \ref{pair} and therefore or Theorem~\ref{equals}.

\subsection{The Calabi-Yau Cap}\label{CYcapsec}

First of all let us notice that we can obtain  $A_d(0|-1,0)_\eta$ from $A_d(0|0,-1)_\eta$  by simply interchanging the roles of $s_1$ and $s_2$. 

\begin{theorem}
\label{CYcap}
 Let $d$ be a positive integer, and $\eta=(\eta_1,\ldots,\eta_{\ell(\eta)})$ a partition of $d$.

 The degree $d$ Calabi-Yau invariants are :
\begin{eqnarray}\begin{vimp}\displaystyle{
A_d(0|0,-1)_{\eta}=(-1)^{d-\ell(\eta)}\frac{\left(2\sin\left(\frac{u}{2}\right)\right)^d}{(s_1)^{\ell(\eta)} \mathfrak{z}{(\eta)}\prod 2\sin\left(\frac{\eta_iu}{2}\right)}.}\end{vimp}\nonumber\end{eqnarray}
\end{theorem}
\textbf{Note:} In \cite{r:adm}, the above formula is computed via localization on moduli spaces of (connected) admissible covers in degree $1,2,3$. The result is obtained by finding  relations between the  Calabi-Yau cap invariants and generating functions for simple Hurwitz numbers. There are two types of obstructions that arise in degree $\geq 4$:
\begin{enumerate}
  \item Fixed loci inside moduli spaces of connected admissible covers are in principle easily described as finite products and quotients of moduli spaces of connected admissible covers, but the combinatorial complexity grows fast.
    \item Generating functions for simple Hurwitz numbers are not readily available beyond degree $3$.
 \end{enumerate}
To circumvent the first problem we interpret the fixed loci in the localization as simpler products of  disconnected admissible cover spaces. Then all possible Calabi-Yau invariants, not only the fully ramified ones, appear in the recursions. There is one subtelty to be aware of: Calabi-Yau cap invariants are defined as intersection numbers on moduli spaces of admissible covers of a parametrized $\proj$, whereas the fixed loci are in terms of admissible covers of unparamterized projective lines. Another localization computation, with an appropriate choice of linearizations for the bundles, gives an expression for the invariants in terms of the unparametrized $\proj$ admissible covers.\\
To deal with the lack of explicit generating functions for general simple Hurwitz numbers, we notice that the recursive relation that we need to prove is in fact determined by a virtual localization computation on moduli spaces of stable maps. This is yet more evidence of how intimately related this theory and Gromov-Witten theory are.

\subsection{Proof of Theorem \ref{CYcap}}
\subsubsection*{The Connected Calabi-Yau Cap invariants}
We prove the following formula for the connected invariants:
\begin{eqnarray}
\begin{vimp}\displaystyle{
\begin{array}{ccc}
A^{\circ}_d(0|0,-1)_{\eta} & = & \left\{\begin{array}{cc}
                              \displaystyle{\frac{(-1)^{d-1}}{s_1}\frac{1}{d}\frac{\left(2\sin\left(\frac{u}{2}\right)\right)^d}{2\sin\left(\frac{du}{2}\right)}} & \mbox{for} \ \ \eta=(d) \\ & \\
0 & \mbox{otherwise}.
                                 \end{array}\right.
\end{array}}
\end{vimp}\label{CYinv}
\end{eqnarray}

Theorem~\ref{CYcap} follows from  (\ref{CYinv}) via exponentiation.

\vspace{0.5cm}
 The vanishing of the connected invariants for  all partitions  but  $(d)$ is a dimension count. By definition (\ref{adminv}) and formula (\ref{ordchern}), the genus $h$ contribution to the connected Calabi-Yau invariants is:

\begin{eqnarray}
{A^\circ}^h_d(0|0,-1)_{\eta}&=&\int_{\overline{Adm}^\circ_{h\stackrel{d}{\rightarrow}\proj,(\eta\infty)}}\hspace{-1cm} c_{2h+d-1}^{eq} (-R^\bullet\pi_\ast f^\ast(\mathcal{O}_{\proj}\oplus \mathcal{O}_{\proj}(-1)))\nonumber\\
 &=& \sum_{b_1,b_2}s_1^{r_1-b_1}s_2^{r_2-b_2}\int_{\overline{Adm}^\circ_{h\stackrel{d}{\rightarrow}\proj,(\eta\infty)}}\hspace{-1cm} c_{b_1}(\mathbb{E}^\ast)c_{b_2}(R^1\pi_\ast f^\ast(\mathcal{O}_{\proj}(-1))),\label{intgen}\end{eqnarray}
where:
\begin{itemize}
         \item $b_1+b_2 = dim (\overline{Adm}_{h\stackrel{d}{\rightarrow}\proj,(\eta\infty)})= 2h + d +\ell(\eta)-2$;
           \item  $r_1=h-1$ is the virtual rank of the virtual bundle $-R^\bullet\pi_\ast f^\ast(\mathcal{O}_{\proj})$;
             \item $r_2=h+d-1$ is the virtual rank of the virtual bundle $-R^\bullet\pi_\ast f^\ast(\mathcal{O}_{\proj}(-1))$.
\end{itemize}  

Since 
$$-R^\bullet\pi_\ast f^\ast(\mathcal{O}_{\proj}(-1))=R^1\pi_\ast f^\ast(\mathcal{O}_{\proj}(-1))$$
is in fact a vector bundle of rank $h+d-1$, we also have the constraint 
$$b_1+b_2 \leq 2h+d-1.$$ 
The only possibly nonvanishing integrals occur when $\ell(\eta)=1$, i.e. when $\eta =(d)$. The indices $b_1$ and $b_2$  are forced to be, respectively, $h$ and $h+d-1$.
\vspace{0.5cm}

\textbf{Remark:} The full ramification condition forces all covers to be connected; the fully ramified connected and disconnected invariants coincide, thus allowing us to drop the superscript ``$\circ$''.
\vspace{0.5cm} 

Finally, our task is to prove:
$$ \frac{1}{s_1}\sum_{h=0}^{\infty} u^{2h+d-1}\int_{\overline{Adm}_{h\stackrel{d}{\rightarrow}\proj,((d)\infty)}}\hspace{-1cm} c_{h}(\mathbb{E}^\ast)c_{h+d-1}( R^1\pi_\ast f^\ast(\mathcal{O}_{\proj}(-1)))=\frac{(-1)^{d-1}}{s_1}\frac{1}{d}\frac{\left(2\sin\left(\frac{u}{2}\right)\right)^d}{2\sin\left(\frac{du}{2}\right)}.$$

\subsubsection*{Calabi-Yau Cap Invariants: Parametrized to Unparametrized}
We evaluate via localization the Calabi-Yau cap invariant $A^h_d(0|0,-1)_{\eta}$, for a general partition $\eta$.

We linearize the ($\Cstar$ action on the) two bundles as indicated in the following table:
\begin{center}
\begin{tabular}{|l||c|c|}
\hline
weight : & over $0$ & over $\infty$ \\
\hline
\hline
 ${\mathcal{O}_{\proj}}(-1)$ & 0 & 1 \\
\hline
${\mathcal{O}_{\proj}}$  & 0  & 0  \\
\hline
\end{tabular}
\end{center}

There are a priori many fixed loci in the localization computation. However it is possible to rule out a vast majority of them using either dimension counts or linearization considerations ( see \cite{r:adm} or \cite{bp:tqft} for a discussion of these standard localization ``tricks'').

Eventually, the only possibly contributing fixed loci are those whose general element consists 
of $\ell(\eta)$ spheres $S_i$, mapping to the main $\proj$ with degree $\eta_i$,  all fully ramified over $0$ and $\infty$. A genus $0$ twig sprouts from the point $\infty$ on the main $\proj$, covered by $\ell(\eta)$ curves $C_i$ of genus $h_i$. The curve $C_i$ is attached to $S_i$ at a fully ramified point. The $h_i$'s are such that
\begin{eqnarray}
h_1 + \ldots + h_{\ell(\eta)} = h + \ell(\eta) -1.
\end{eqnarray}  
Finally, if we denote by $F_{\eta,h}$ the disjoint union of all such fixed loci as the $h_i$'s vary, and by $N$ the normal bundle to such fixed loci, we obtain from localization:
\begin{eqnarray}
A^h_d(0|0,-1)_{\eta}=\int_{F_{\eta,h}}\frac{e^{eq} (-R^\bullet\pi_\ast f^\ast(\mathcal{O}_{\proj}\oplus \mathcal{O}_{\proj}(-1)))\mid_{F_{h,\eta}}}{e^{eq}(N)}.
\end{eqnarray}

\subsubsection*{Recursion via Localization on Admissible Covers}

We now suppose $d>1$ and consider the following auxiliary integral, computed on the space of connected admissible covers:
\begin{eqnarray}
I^h=\int_{\overline{Adm}^\circ_{h\stackrel{d}{\rightarrow}\proj}}\hspace{0cm} e^{eq} (-R^\bullet\pi_\ast f^\ast(\mathcal{O}_{\proj}\oplus \mathcal{O}_{\proj}(-1))).
\label{aux}
\end{eqnarray}

Elementary dimension reasons give us the vanishing of (\ref{aux}): we are integrating a class whose highest non-equivariant factor has codimension $(2h+d-1)$ on a space of dimension $2h+2d-2$.

On the other hand, if we evaluate the integral via localization we  get a relation among Calabi-Yau cap invariants. We let a one dimensional torus act naturally on the moduli space and denote the equivariant parameter  $s$. We choose to linearize the two bundles with weights:

\begin{center}
\begin{tabular}{|l||c|c|}
\hline
weight : & over $0$ & over $\infty$ \\
\hline
\hline
 ${\mathcal{O}_{\proj}}(-1)$ & -1 & 0 \\
\hline
${\mathcal{O}_{\proj}}$  & 1  & 1  \\ \hline
\end{tabular}
\end{center}

The possibly contributing fixed loci $E_{\eta,h_0,h_\infty}$ are represented by connected localization graphs such that any vertex over $\infty$ has valence $1$ ( \cite{r:adm} ). They can be indexed by triples $(\eta, h_0,h_\infty)$, where 
\begin{itemize}
  \item $\eta=(d_1,\ldots, d_{\ell(\eta)})$ is a partition of $d$ representing the configuration of the spheres over the main $\proj$;
    \item $h_0$ is the genus of the curve lying over $0$;
      \item $h_\infty$ is the genus of the curve lying over $\infty$ (considered as a disconnected curve);
        \item $h_0+h_\infty= h - \ell(\eta) +1$.
\end{itemize}

We recognize that a general element in the fixed locus $E_{\eta,h_0,h_\infty}$ is obtained by gluing together an element in the fixed locus $F_{\eta,h_\infty}$ with a  connected admissible covers of a genus $0$ curve, with a special point of ramification $\eta$. Keeping in account the  stacky contribution from the gluing, then our integral $I$ on $E_{\eta,h_0,h_\infty}$ reduces to:
$$I^h_{\eta,h_0,h_\infty}=
\mathfrak{z}(\eta)\int_{\overline{Adm}^\circ_{h_0\stackrel{d}{\rightarrow}\proj,\eta}\times F_{\eta,h_\infty}} \hspace{-0.5cm}\frac{e^{eq} (-R^\bullet\pi_\ast f^\ast(\mathcal{O}_{\proj}\oplus \mathcal{O}_{\proj}(-1)))\mid_{\left( \overline{Adm}^\circ_{h\stackrel{d}{\rightarrow}\proj,\eta}\times F_{\eta,h_\infty}\right)}}{e^{eq}(N)}=
$$
$$=\mathfrak{z}(\eta)s^{2\ell(\eta)} A_d^{h_\infty}(0|0,-1)_\eta\int_{\overline{Adm}^\circ_{h_0\stackrel{d}{\rightarrow}\proj,\eta}}\frac{c_{h_0}(\mathbb{E^\ast}\otimes \mathbb{C}_1)c_{h_0}(\mathbb{E^\ast}\otimes \mathbb{C}_{-1})}{s(s-\psi_\eta)},$$
where $\mathbb{C}_a$ denotes a trivial line bundle where the torus acts on the fibers with weight $a$.

After expanding the above expression and simplifying using Mumford's relation (\cite{m:taegotmsoc})
\begin{eqnarray}
c(\mathbb{E})c(\mathbb{E}^\ast)=1,
\label{mum}
\end{eqnarray}
we obtain
\begin{eqnarray}I^h_{\eta,h_0,h_\infty}=
\mathfrak{z}(\eta)s^{\ell(\eta)+2-d} A_d^{h_\infty}(0|0,-1)_\eta\int_{\overline{Adm}^\circ_{h_0\stackrel{d}{\rightarrow}\proj,\eta}}(-1)^{h_0}\psi_\eta^{2h_0 +d +\ell(\eta)-4} = \nonumber \\
=\mathfrak{z}(\eta)s^{\ell(\eta)+2-d} A_d^{h_\infty}(0|0,-1)_\eta\frac{(-1)^{h_0}H_d^{h_0}({\eta})}{(2h_0 +d +\ell(\eta)-2)!}.
\end{eqnarray}
The quantity $H_d^{h_0}({\eta})$ is a simple Hurwitz number, as defined in section $\ref{dijkgraaf}$.

The evaluation of the integral $I$ is obtained by adding up the contributions coming from all fixed loci $E_{\eta,h_0,h_\infty}$:
\begin{eqnarray}
0=I^h= \sum_{\eta \vdash d}\sum_{h_0+h_\infty = h -\ell(\eta)+1} I^h_{\eta,h_0,h_\infty}.
\label{eval}
\end{eqnarray}

Formula (\ref{eval}) holds for all genera $h$. All such formulas can be expressed in a very compact form in the language of generating functions. Define:
\begin{itemize}
  \item $\displaystyle{\mathcal{H}_{d,\eta(u)}:=\sum\frac{(-1)^{h}H_d^{h}({\eta})}{(2h +d +\ell(\eta)-2)!} u^{(2h +d +\ell(\eta)-2)}}.$
\end{itemize}

Then formulas (\ref{eval}), for all genera $h$, are encoded in the relation:
\begin{eqnarray}
0=\sum_{\eta \vdash d} \mathfrak{z}(\eta)s^{\ell(\eta)+2-d} A_d(0|0,-1)_\eta(u) \mathcal{H}_{d,\eta}(u).
\label{fundrel}
\end{eqnarray} 

Relation (\ref{fundrel}) determines $A_d(0|0,-1)_{(d)}$ in terms of generating functions for simple Hurwitz numbers and of the invariants $A_d(0|0,-1)_\eta$, for $\ell(\eta)\geq 2$, which can be inductively determined via exponentiation if we assume  the theory up to degree $d-1$. The theory has been explicitly computed up to degree $3$ in~\cite{r:adm}, hence the induction can start.

To prove Theorem~\ref{CYcap} it therefore suffices to show that relation (\ref{fundrel}) holds for the conjectured values of the Calabi-Yau invariants. After substituting and simplifying, this amounts to proving:
\begin{eqnarray}
0=\sum_{\eta\vdash d}(-1)^{\ell(\eta)}\frac{\mathcal{H}_{d,\eta}(u)}{\displaystyle{\prod_{\eta_i\in\eta} 2\sin\left(\frac{\eta_iu}{2}\right)}}.
\label{relfin}
\end{eqnarray}

\subsubsection*{Virtual Localization on Stable Maps} 

Relation (\ref{relfin}) is the result of explicitly evaluating via virtual localization the following auxiliary integrals:
\begin{eqnarray}
J^h=\int_{\overline{M}_{h}(\proj,d)}\hspace{0cm} e^{eq} (-R^\bullet\pi_\ast f^\ast(\mathcal{O}_{\proj}\oplus \mathcal{O}_{\proj}(-1))).
\label{aux2}
\end{eqnarray}  

Again dimension reasons grant us the vanishing of this integral. We proceed to linearize the bundles as shown in the following table:
\begin{center}
\begin{tabular}{|l||c|c|}
\hline
weight : & over $0$ & over $\infty$ \\
\hline
\hline
 ${\mathcal{O}_{\proj}}(-1)$ & -1 & 0 \\
\hline
${\mathcal{O}_{\proj}}$  & 1  & 1  \\ \hline
\end{tabular}
\end{center}

The analysis of the possibly contributing fixed loci is parallel to the previous section. The contribution by the fixed locus $E_{\eta,h_0,h_\infty}$ is (see \cite{mirror}, chapter 27, for a clear and detailed explanation of how to compute these terms, or~\cite{bp:tlgwtoc}, proof of Theorem $5.1$ for an extremely similar computation):
$$
\sum_{h_1+\ldots+h_{\ell(\eta)}= h_\infty+\ell(\eta)-1} J_{\eta,h_0,h_1,\ldots,h_{\ell(\eta)}},
$$
with
\begin{eqnarray}
 J_{\eta,h_0,h_1,\ldots,h_{\ell(\eta)}}= 
\frac{1}{\mathfrak{z}(\eta)}
\int_{\overline{M}_{h_0,\ell(\eta)}}
\frac{c_{h_0}(\mathbb{E}^\ast\otimes \mathbb{C}_1)c_{h_0}(\mathbb{E}^\ast\otimes \mathbb{C}_1)c_{h_0}(\mathbb{E}^\ast\otimes \mathbb{C}_{-1})}
{\prod \left(\frac{1}{\eta_i}-\psi_i\right)} \nonumber \\
\prod_{i=1}^{\ell(\eta)}\frac{\eta_i^{\eta_i}}{\eta_i!}
\int_{\overline{M}_{h_i,1}}
\frac{c_{h_i}(\mathbb{E}^\ast\otimes \mathbb{C}_1)c_{h_i}(\mathbb{E}^\ast\otimes \mathbb{C}_1)c_{h_i}(\mathbb{E}^\ast)}
{-\frac{1}{\eta_i}-\psi_1}.\label{loc}
\end{eqnarray}

After simplifying via Mumford's relation and rearranging things, formula (\ref{loc}) becomes:

\begin{eqnarray}
\frac{(-1)^{h_0}}{Aut(\eta)}
\prod_{i=1}^{\ell(\eta)}\frac{\eta_i^{\eta_i}}{\eta_i!}
\int_{\overline{M}_{h_0,\ell(\eta)}}
\frac{1-\lambda_1 + \ldots \pm \lambda_{h_0}}
{\prod (1-{\eta_i}\psi_i)} 
\prod_{i=1}^{\ell(\eta)}-\eta_i^{2h_i-1}
\int_{\overline{M}_{h_i,1}} \lambda_{h_i}\psi_1^{2h_i-2}.
\label{score}
\end{eqnarray}

We recognize in formula (\ref{score}) two famous results in the field:
\begin{description}
\item[ELSV formula:] in \cite{elsv:hnaiomsoc} and \cite{vg:hnavl}, this formula estabilishes the 
connection between Hurwitz numbers and Hodge integrals:
\begin{eqnarray}
H^h_d(\eta)=
\frac{(2h+d+\ell(\eta)-2)!}{Aut(\eta)}
\prod_{i=1}^{\ell(\eta)}\frac{\eta_i^{\eta_i}}{\eta_i!}
\int_{\overline{M}_{h,\ell(\eta)}}
\frac{1-\lambda_1 + \ldots \pm \lambda_{h}}
{\prod (1-{\eta_i}\psi_i)} 
\label{ELSV}
\end{eqnarray}
\item[Faber and Pandharipande's formula:] in \cite{fp:hiagwt}, the following class of integrals is
computed and expressed in generating function form:
\begin{eqnarray}
\mathcal{L}(u):= \sum u^{2h-1}\int_{\overline{M}_{h,1}} \lambda_{h}\psi_1^{2h-2} 
= \frac{1}{2 \sin\left(\frac{u}{2}\right)}.
\label{fapa}
\end{eqnarray}
\end{description}
Now it is a matter of careful bookkeeping to translate all of the above information in the language of generating functions.  
Doing so concludes the proof of Theorem \ref{CYcap} by estabising the truth of relation
(\ref{relfin}):
\begin{eqnarray}
0=\sum_{h\in\mathbb{Z}} J^hu^{2h+2d-2} = \sum_{\eta \vdash d} (-1)^{\ell(\eta)}
\mathcal{H}_{d,\eta}(u)\prod_{\eta_i\in\eta} \mathcal{L}(\eta_iu) =  \nonumber \\
=
\sum_{\eta\vdash d}(-1)^{\ell(\eta)}\frac{\mathcal{H}_{d,\eta}(u)}{\displaystyle{\prod_{\eta_i\in\eta} 2\sin\left(\frac{\eta_iu}{2}\right)}} \ \ .
\end{eqnarray}

\section{A Specialization of the Theory}\label{special}
We now  discuss a specialization of the theory, obtained by embedding 
a one-dimensional torus inside the two-dimensional torus $T$,
 and considering the theory as depending from one equivariant
 parameter instead of two. 

We specialize to the anti-diagonal action, and notice that the coefficients for
 the product simplify dramatically.
 It is  possible to obtain nice closed formulas for our theory,
 and to view our TQFT as a one parameter deformation of the classical TQFT
 of Hurwitz numbers studied by Dijkgraaf and Witten in \cite{dw:tgtagc} and Freed and Quinn in \cite{fq:cst}. 
Our formulas show connections to the representation theory of the symmetric group $S_d$. The relevant representation theoretic quantities  are introduced in  Appendix~\ref{app1}.

\subsection{The Anti-diagonal Action}

Let $\Cstar$ be embedded in the two-dimensional torus $T$ via the map
$$\alpha  \mapsto \left( \alpha, \frac{1}{\alpha}\right).$$

$\Cstar$ acts on $N$ by composing this embedding with the natural action of $T$ constructed in page \pageref{toract}. If we let
$$H^\ast_{\Cstar}(pt)= \mathbb{C}[s],$$
then  the one parameter theory obtained with this action corresponds to setting
$$s=s_1=-s_2.$$

\subsubsection{The Q-dimension of an Irreducible Representation}
Let $\rho$ be an irreducible representation of the symmetric group on $d$ letters $S_d$. Classically, a Young diagram, and hence a partition of $d$, can be canonically associated to  $\rho$ (see Appendix~\ref{app1}).

We now define the $Q$-dimension of the representation $\rho$ to be:
\begin{eqnarray}
\frac{dim_Q\rho}{d!}:= \prod_{\Box\in\rho}\frac{1-Q}{1-Q^{h(\Box)}}=\prod_{\Box\in\rho}\frac{1}{1+Q+\dots+Q^{h(\Box)-1}}
\label{Qdim}\end{eqnarray}

As a consequence of the classical hooklenght formula (\ref{hooklength}), formula (\ref{Qdim}) specializes to the ordinary dimension of $\rho$ when setting $Q=1$.

\subsubsection{The Level $(0,0)$ TQFT}
The main result is that the level $(0,0)$ TQFT completely collapses to the Dijkgraaf TQFT  $\mathcal{D}$. In particular,
we have  explicit formulas for the semisimple basis of the Frobenius algebra. The basis vectors are indexed by irreducible representations of the symmetric group $S_d$.
\begin{lemma}[essentially Bryan-Pandharipande]
For the anti-diagonal action, the level $(0,0)$ series have no nonzero terms of positive degree in $u$.
\end{lemma}
\textsc{Proof:} endow $\mathbb{C}$ with the  $\Cstar$ action
$$\alpha \cdot z= \alpha^nz.$$
This corresponds to considering $\mathbb{C}$ as an equivariant line bundle over a point, whose first equivariant Chern class is $ns$. We  denote such equivariant line bundle by $\mathbb{C}_{ns}$.

The level $(0,0)$ partition functions are, up to some pure weight factor, constructed from integrals of the form:
$$\int_{\overline{Adm}_{h\stackrel{d}{\longrightarrow}X,(\eta^1x_1,\dots,\eta^rx_r)}}\hspace{-1.5cm}e^{eq}(\mathbb{E}^\ast\otimes\mathbb{C}_s)e^{eq}(\mathbb{E}^\ast\otimes\mathbb{C}_{-s})=$$

$$\int_{\overline{Adm}_{h\stackrel{d}{\longrightarrow}X,(\eta^1x_1,\dots,\eta^rx_r)}}\hspace{-1.5cm}(-1)^he^{eq}((\mathbb{E}^\ast\oplus\mathbb{E})\otimes\mathbb{C}_s).$$
Equivariant Chern classes of a bundle also are products of ordinary Chern classes times the appropriate factor of $s$. But by Mumford's relation (\ref{mum}), all Chern classes (but the $0$-th) of the bundle $\mathbb{E}^\ast\oplus\mathbb{E}$ vanish.
Hence the only possibly nonvanishing integrals  occur when the dimension of the moduli space is $0$, which then constitutes the degree $0$ term in our generating functions.

We have therefore already essentially  produced  a semisimple basis for the corresponding Frobenius algebra in page \pageref{ssb}. All we need to do is to adjust for the equivariant parameter:
\begin{description}
  \item[semisimple basis:] let $\rho$ be an irreducible representation of the symmetric group $S_d$, $\chi_{\rho}$ its character function;  a semisimple basis for the level $(0,0)$ TQFT is given by the vectors
\begin{eqnarray}
\mbox{e}_{\rho}= \frac{dim\rho}{d!}\sum_{\eta \vdash d} (s)^{\ell(\eta)-d}\chi_{\rho}(\eta) \mbox{e}_\eta.
\end{eqnarray}\end{description}
\subsubsection{The Structure of the Weighted TQFT}
\begin{theorem}\label{antid}
The closed partition functions for the weighted TQFT are given by the following closed formulas:

\begin{eqnarray}
\begin{imp}\small{
A_d(g|k_1,k_2)=(-1)^{a} s^{b}\sum_{\rho}\left(\frac{d!}{dim\rho}\right)^{2g-2} \left(\frac{dim\rho}{dim_Q\rho}\right)^{k_1+k_2} Q^{n(\rho)k_1+n(\rho')k_2},
}\end{imp} \nonumber \\
 \nonumber  \\ 
\end{eqnarray}
where:
\begin{itemize}
\item $a= d(g-1-k_2)$;
\item $b=d(2g-2-k_1-k_2)$;
  \item the variable $Q= e^{iu}$;
  \item $\rho$ denotes an irreducible representation of the symmetric group $S_d$;
    \item $\rho'$ denotes the dual representation;
      \item the function $n$ is defined in  section~\ref{n}.
\end{itemize}
\end{theorem}

\textbf{Remark.} Notice that by setting $Q=1$, which corresponds to $u=0$, we recover the classical formula (\ref{burn}) counting unramified covers of a genus $g$ topological surface. Thus any  TQFT naturally embedded in our weighted TQFT constitutes a one parameter deformation of the Dijkgraaf TQFT.

\textsc{Proof:} by Fact~\ref{wtqftstructure}, to completely describe the structure  of a semisimple weighted TQFT it suffices to evaluate the following quantities:
\begin{description}
  \item[$\lambda_\rho$:] the e$_\rho$-eigenvalue of the  genus adding operator, or, equivalently, the inverse of the counit evaluated on e$_\rho$;
    \item[$\mu_\rho$:] the e$_\rho$-eigenvalue of the left level-subtracting operator, or, equivalently, the coefficient  of e$_\rho$  in the $(0,-1)$ Calabi-Yau cap;
 \item[$\overline{\mu}_\rho$:] the e$_\rho$-eigenvalue of the right level-subtracting operator, or, equivalently, the coefficient  of e$_\rho$  in the $(-1,0)$ Calabi-Yau cap.
\end{description}

The computation of $\lambda_\rho$ coincides exactly with  Bryan and Pandharipande's computation in \cite{bp:tlgwtoc}. We reproduce it here  for the sake of completeness.
$$
\begin{array}{rcl}
 \lambda_\rho ^{-1}& =& \mathcal{U}\left(\begin{array}{c}\includegraphics[height=1cm]{counit.eps}\\ (0,0)\end{array} \right)(\mbox{e}_\rho)  \nonumber \\ & & \\
                   & =& \displaystyle{\frac{dim\rho}{d!}\sum_\eta (is)^{\ell(\eta)-d}\chi_\rho (\eta) A_d(0|0,0)_\eta} \nonumber \\ & & \\
                   & = &\displaystyle{\frac{dim\rho}{d!} (is)^{\ell(1^d)- d}\chi_\rho(1^d) \frac{1}{d! (-s^2)^d}} \nonumber \\ & & \\
                   & = &\displaystyle{\left(\frac{dim\rho}{d!}\right)^2 (is)^{-2d}} \nonumber
\end{array}
$$
Hence,
\begin{eqnarray}
 \lambda_\rho = \left(\frac{d!}{dim\rho}\right)^2 (is)^{2d} \nonumber
\end{eqnarray}

To compute $\mu_\rho$ and $\overline{\mu}_\rho$ let us first of all observe that the  tensors  associated to the Calabi-Yau caps  in our theory are scalar multiples of the tensors in Bryan and Pandharipande's theory.
$$\mathcal{U}(CY cap) = 2^d\sin\left(\frac{u}{2}\right)^d \mathcal{BP}(CY cap)= \frac{(1-Q)^d}{Q^\frac{d}{2}(-i)^d}\mathcal{BP}(CY cap). $$
This observation, together with the formulas  in \cite{bp:tlgwtoc}, page 36, implies:
$$\mu_\rho= s^d\frac{d!}{dim \rho}(1-Q)^d s_\rho(Q),$$

$$\overline{\mu}_\rho= (-s)^d\frac{d!}{dim \rho}(1-Q)^ds_{\rho'}(Q) ,$$
where $s_\rho$ denotes the Schur function of the representation $\rho$, and  is defined to be (\cite{m:sfahp}):
$$ s_\rho(Q):= Q^{n(\rho)}\prod_{\Box \in \rho} \frac{1}{1-Q^{h(\Box)}}.$$
Plugging this in, we obtain:
$$
\begin{array}{rcl}
\mu_\rho & = & \displaystyle{s^d\left(\frac{d!}{dim \rho}\right)(1-Q)^d Q^{n(\rho)}\prod_{\Box \in \rho} \frac{1}{1-Q^{h(\Box)}}}\\ & & \\
         & = & \displaystyle{s^d\left(\frac{d!}{dim \rho}\right)Q^{n(\rho)}\prod_{\Box \in \rho} \frac{1-Q}{1-Q^{h(\Box)}}}\\ & & \\
         & = & \displaystyle{s^d\left(\frac{dim_Q\rho}{dim \rho}\right)Q^{n(\rho)}}\\
\end{array}
$$
and
$$
\begin{array}{rcl}
\overline{\mu}_\rho & = & \displaystyle{(-s)^d\left(\frac{d!}{dim \rho}\right)(1-Q)^d Q^{n(\rho')}\prod_{\Box \in \rho'} \frac{1}{1-Q^{h(\Box)}}}\\ & & \\
         & = & \displaystyle{s^d\left(\frac{d!}{dim \rho}\right)Q^{n(\rho')}\prod_{\Box \in \rho'} \frac{1-Q}{1-Q^{h(\Box)}}}\\ & & \\
         & = & \displaystyle{s^d\left(\frac{dim_Q\rho}{dim \rho}\right)Q^{n(\rho')}}\\
\end{array}
$$
Theorem~\ref{antid} is finally obtained by using these coefficients in the formula given by Fact~\ref{wtqftstructure}.

\appendix
\section{Combinatorics and Representation Theory} \label{app1}
\subsection{Partitions of an Integer\label{n}}

A partition $\eta$ of an integer $d$ is a finite sequence of positive integers adding up to $d$:
$$\eta= (\eta^1,    \eta^2,    \dots,\eta^r),$$
with
$$\eta^1\geq\eta^2\geq\dots\geq\eta^r$$
and
$$|\eta|=\sum_i \eta^i=d.$$
We use the notation $\eta\vdash d$ to indicate that $\eta$ is a partition of $d$.

The number $r$ of nonzero integers is called the \textit{length} of the partition $\eta$, and denoted  $\ell(\eta)$.

It is also convenient to have a notation that groups all equal parts together. By
$$((\eta^1)^{m_1}\dots(\eta^k)^{m_k})$$
we denote the partition
$$\eta=(\underbrace{\eta^1,\dots,\eta^1}_{m_1\ times},\underbrace{\eta^2,\dots,\eta^2}_{m_2\ times},\dots,\underbrace{\eta^k,\dots,\eta^k}_{m_k\ times}).$$

A partition $\eta$ of the integer $d$ can be canonically represented by a Young diagram. A \textit{Young diagram} is a left justified array of boxes, such that the lenght of the rows does not increase as you go down the diagram. To a partition $\eta$ we associate the Young diagram whose $i$-th row is composed of $\eta^i$ boxes. For example:
$$(3,2,2,1,1)= (3^1 2^2 1^2) = \parbox[c]{0.10\textwidth}{\includegraphics[width=1cm]{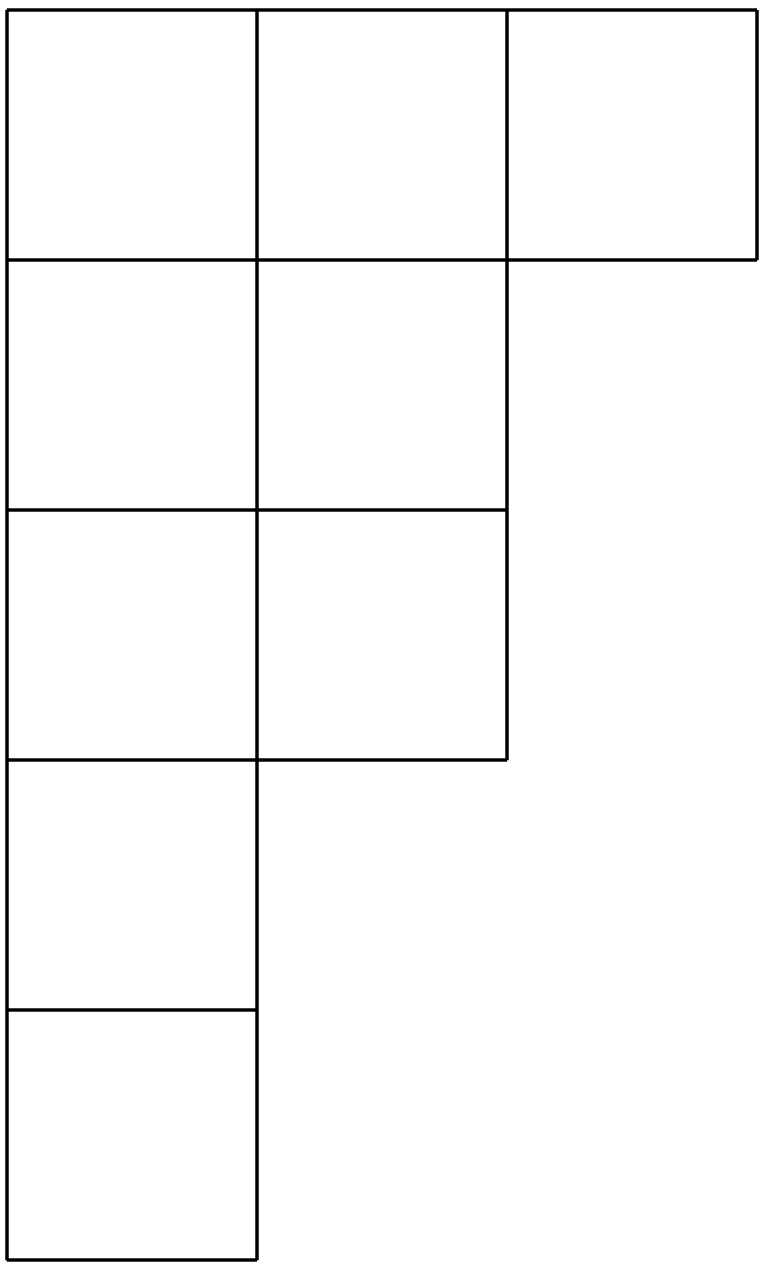}}.$$

The conjugate partition $\eta'$ is the partition associated to the reflection along the main diagonal of the Young diagram associated to $\eta$. In our example:
$$(3,2,2,1,1)'= (5,3,1)= \parbox[c]{0.15\textwidth}{\includegraphics[height=1cm]{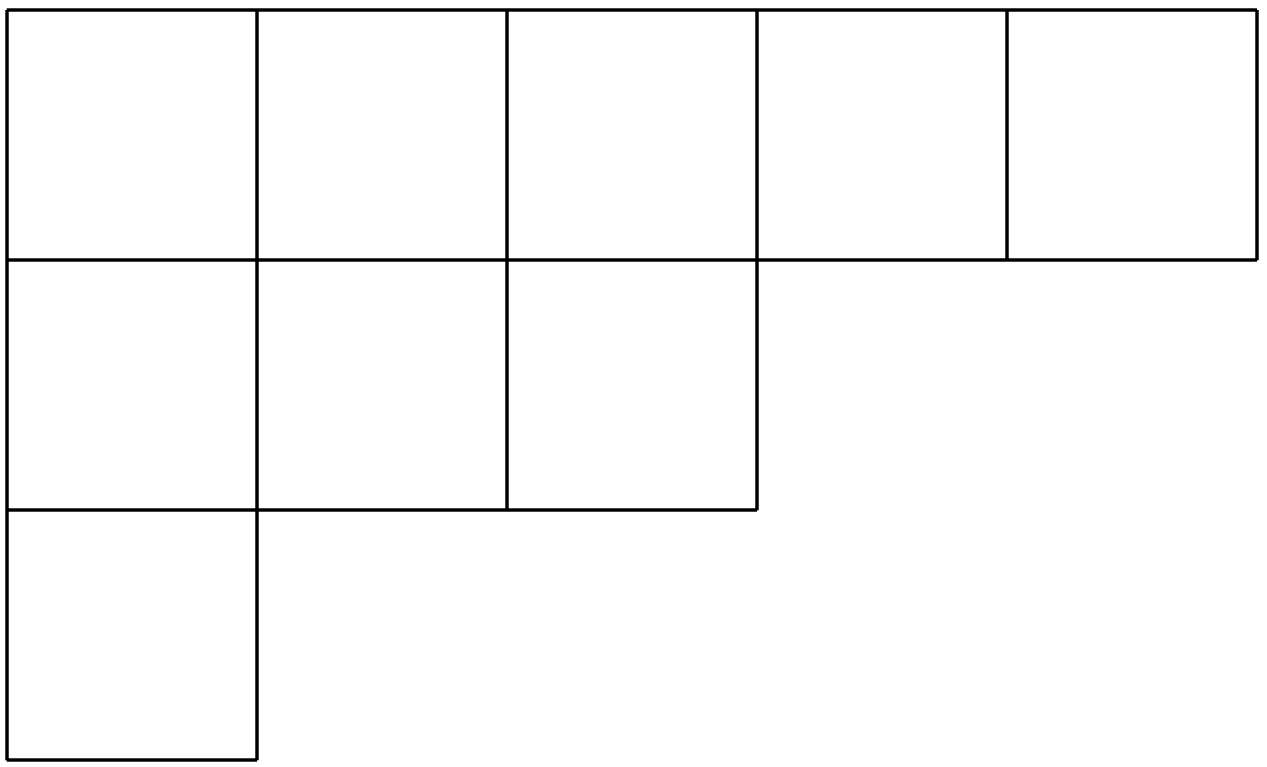}}.$$

Here are some functions associated to Young diagrams that are useful for our purposes:
\begin{description}
  \item[content:] for a given box $\Box$ in position $(i,j)$ in a Young diagram, we define the content $c(\Box)=i-j$. Boxes on the main diagonal  have content $0$, boxes above the diagonal  have positive content that measures exacly how many diagonals over the main one they lie on, and so on.
The \textit{total content} $c(\eta)$ is defined to be the sum of the content of all boxes in the diagram. If a diagram is symmetric with respect to the main diagonal, then its total content is $0$. In some sense, the total content measures the asymmetry of a Young diagram.

\item[$n$-function:] given a Young diagram, we define the $n$-function as follows: 
number all boxes in the first row with $0$'s, all boxes in the second row with $1$'s,
 all boxes in the third row with $2$'s and so on. 
Then add all of these numbers to obtain $n(\eta)$. 

\end{description}

The following obvious formula connects these two quantities:
$$c(\eta)= n(\eta')-n(\eta).$$
\begin{description}
\item[hooklength:] for a given box $\Box$ in a Young diagram, the \textit{hooklength} $h(\Box)$ is the length (as in number of boxes) of the hook that has the given box as its north-west corner, as shown in Figure~\ref{hook}.

\begin{figure}[htbp]
\begin{center}
\includegraphics[width= 3cm]{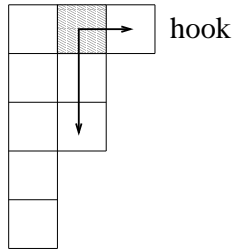}
\caption{The hook associated to the shaded box.}\label{hook}
\end{center}
\end{figure}
\end{description}

The total hooklength $h(\eta)$ is the sum of the hooklengths over all boxes in the diagram.

\subsection{The Hooklength Formula} Let $\rho$ be an irreducible representation of $S_d$  represented by the Young diagram associated to the partition $\eta$. Then the following formula holds:
\begin{eqnarray}\mbox{dim}\rho= \frac{\displaystyle{d!}}{\displaystyle{\prod_{\Box\in \eta}h(\Box)}}. \label{hooklength}\end{eqnarray}

\addcontentsline{toc}{section}{Bibliography}
\bibliographystyle{alpha}
\bibliography{biblio}

\end{document}